\documentclass{article}

\usepackage[english]{babel}

\usepackage[letterpaper,top=2cm,bottom=2cm,left=3cm,right=3cm,marginparwidth=1.75cm]{geometry}

\usepackage{fullpage}
\usepackage{amsfonts}
\usepackage{graphicx}
\usepackage{amsmath}
\usepackage{amsthm}
\usepackage{amssymb}
\usepackage{bm}
\usepackage[colorlinks=true, allcolors=blue]{hyperref}
\usepackage[english]{babel}
\usepackage[letterpaper,top=2cm,bottom=2cm,left=3cm,right=3cm,marginparwidth=1.75cm]{geometry}
\usepackage{xcolor}
\usepackage{mathabx}
\usepackage{tikz}
\usetikzlibrary{patterns, arrows.meta, positioning, decorations.pathreplacing, shapes, fit, backgrounds}
\usepackage{enumitem}

\newtheorem{example}{Example}

\newcommand{\bx}{\bm{x}}
\newcommand{\bbeta}{\bm{\beta}}
\newcommand{\mX}{\bm{X}}
\newcommand{\bv}{\bm{v}}
\newcommand{\bu}{\bm{u}}
\newcommand{\mZ}{\bm{Z}}
\newcommand{\R}{\mathbb{R}}

\title{High-Dimensional Statistics: Reflections on Progress and Open Problems}
\author{Arian Maleki$^{1}$,~
Subhabrata Sen$^{2}$,~ Sivaraman Balakrishnan$^{3}$,~ Verena Zuber$^{4}$,\\
Chao Gao$^{5}$,
~ Rishabh Dudeja$^{6}$, Christos Thrampoulidis$^{7}$,
~ Anru Zhang$^{8}$,~ Weijie Su$^{9}$,\\~ Jason M. Klusowski$^{10}$,~ Po-Ling Loh$^{11}$,~Ali Shojaie$^{12}$\footnote{The order of authors has been randomized.}}

\begin{document}
\maketitle

\begin{center}
\small
$^{1}$ Department of Statistics, Columbia University\\
$^{2}$ Department of Statistics, Harvard University\\
$^{3}$ Department of Statistics \& Data Science and Machine Learning Department, Carnegie Mellon University\\
$^{4}$ Department of Epidemiology and Biostatistics, School of Public Health, Imperial College London\\
$^{5}$ Department of Statistics, University of Chicago\\
$^{6}$ Department of Statistics, University of Wisconsin--Madison\\
$^{7}$ Department of Electrical and Computer Engineering, University of British Columbia \\
$^{8}$ Department of Biostatistics \& Bioinformatics and Department of Computer Science, Duke University\\
$^{9}$ Department of Statistics and Data Science, Wharton School, University of Pennsylvania\\
$^{10}$ Department of Operations Research and Financial Engineering, Princeton University\\
$^{11}$ Statistical Laboratory, Department of Pure Mathematics and Mathematical Statistics, University of Cambridge\\
$^{12}$ Department of Biostatistics and Department of Statistics, University of Washington, Seattle
\end{center}

\begin{abstract}

Over the past two decades, the field of high-dimensional statistics has experienced substantial progress, driven largely by technological advances that have dramatically reduced the cost and effort for data collection and storage across a broad range of domains, including biology, medicine, astronomy, and the social and environmental sciences. Modern datasets are increasingly complex, often exhibiting rich dependency, heterogeneity, and other features that challenge traditional statistical methods. In response, high-dimensional statistics has evolved to address more sophisticated estimation and inference problems. This evolution has, in turn, fostered deep connections with and contributions to a wide range of research areas, including optimization, concentration of measure, random matrix theory, information theory, and theoretical computer science. Given the rapid pace of recent developments in high-dimensional statistics, our goal is to synthesize representative advances, highlight common themes and open problems, and point to important works that offer entry points into the field. 
\end{abstract}

\begin{sloppypar}

\section{Introduction}
For much of the 19$^{\text{th}}$ and early 20$^{\text{th}}$ centuries, statistical theory was developed with low-dimensional problems in mind---settings in which the number of features was small relative to the number of observations, giving rise to classical assumptions involving fixed dimensions and abundant data. Under such assumptions, statisticians derived elegant, foundational theory for likelihood-based methods such as maximum likelihood estimation and likelihood ratio tests, and established their statistical optimality. These ideas subsequently permeated into data analytic pipelines across the natural and social sciences.

%Consequently, foundational theory was developed for methods such as maximum likelihood estimation and likelihood ratio tests, which became central tools for addressing a wide range of statistical problems.

Rapid technological advances in the late 20$^{\text{th}}$ and early 21$^{\text{st}}$ centuries dramatically reduced the cost and effort required to collect, store, and process data, fundamentally transforming the scale and structure of modern datasets. New technology for tasks such as genome sequencing, imaging, and remote sensing enabled the collection of thousands or even millions of features, thereby challenging statistical methods developed under classical assumptions. Many algorithms and methods derived using fixed-dimensional asymptotics or matrix inversion became unstable or ill-posed in the presence of large numbers of features.

Confronted with these challenges, statisticians systematically developed new theory and computational tools. The collective efforts of researchers across many fields have coalesced into high-dimensional statistics, a vibrant research area focused on datasets in which the number of features is proportional to or exceeds the sample size. Ideas from classical statistics, probability theory, optimization, information theory, computer science, and statistical physics have generated principled methods that remain reliable in high-dimensional settings. We refer the reader to the excellent article \cite{johnstone2009statistical} for a snapshot of the field in its early days.

Developments in high-dimensional statistics have provided a unifying framework for understanding and analyzing complex, large-scale data across diverse domains. In parallel with theoretical advances, the field has introduced a rich set of methodological tools designed to address challenges such as dimensionality exceeding sample size, complex dependence structure, and computational scalability. Prominent examples include sparsity-inducing methods such as the Lasso, related regularization techniques, and modern tools for uncertainty quantification---including debiased estimators and selective inference---which have had a profound impact on downstream applications such as neuroscience, genomics, medical imaging, finance, economics, computational imaging, and signal processing, where extracting interpretable and reliable information from massive datasets is a central challenge.

% High-dimensional statistics has also initiated a conceptual shift in statistics and related fields by transforming how complexity, uncertainty, and inference are understood in modern data analysis. This shift has elevated structure---such as sparsity, low-rankness, and manifold constraints---from a modeling convenience to a central principle that enables meaningful inference in otherwise ill-posed problems. Moreover, high-dimensional statistics has blurred traditional boundaries between statistics, optimization, information theory, and computer science, promoting algorithmic considerations and computational tractability as intrinsic components of statistical reasoning. Beyond methodology, it has reshaped how uncertainty is quantified, how models are validated, and how scientific questions are framed in new data regimes, influencing not only statistics but also the broader landscape of data science and applied scientific inquiry.

High-dimensional statistics has initiated a profound conceptual shift in statistics and related fields, transforming how complexity, uncertainty, and inference are understood in modern data analysis. A central lesson is that dimensionality itself is rarely the fundamental obstacle; rather, meaningful structure, such as sparsity, low rank, or more general notions of low complexity can enable meaningful inference in otherwise ill-posed problems.
This perspective has elevated structure from a modeling convenience to a guiding principle, motivating the widespread use of shrinkage and regularization. In parallel, classical asymptotic reasoning has given way to a non-asymptotic viewpoint centered on finite-sample guarantees and modern proportional asymptotics. High-dimensional statistics has also led to a greater appreciation for the challenges of valid statistical inference, particularly the fact that strong estimation performance does not necessarily imply valid uncertainty quantification. Perhaps most importantly, high-dimensional statistics has blurred traditional boundaries between statistics, optimization, information theory, and computer science, recognizing computation as an intrinsic component of statistical reasoning. For example, to estimate  high-dimensional parameters under structural constraints, statisticians have developed new methods that leverage recent advances in optimization. The statistical performance of these novel methods has been subsequently characterized using tools from information theory and high-dimensional probability. Finally, interactions with theoretical computer science have revealed regimes where optimal statistical performance might be impossible to achieve using computationally tractable learning algorithms. 
Beyond these specific contributions, high-dimensional statistics has reshaped uncertainty quantification,
model validation, and how scientific questions are framed in new data regimes, influencing not only statistics but also the broader landscape of data science and applied scientific inquiry.

% Although high-dimensional statistics has become a standard topic covered in many graduate (or even undergraduate) university curricula, countless open research directions remain. Some questions are driven by the desire to apply a high-dimensional viewpoint to increasingly complex data types, such as those arising from modern AI systems. Other questions stem from the fact that in the modern ``data revolution,'' datasets are so massive that they must be stored/compressed/partitioned/communicated in ways that are parsimonious, yet do not lose important information. If the goal of high-dimensional statistics is to uncover underlying low-dimensional structure in large data sets, the best path forward is not  clear when standard assumptions are no longer met, data come from heterogeneous pools, or structural representations themselves must be learned from the data.

Although high-dimensional statistics is now a standard part of many graduate (or even undergraduate) curricula, the field continues to face a wide range of fundamental challenges. Some arise from the need to extend high-dimensional methods to increasingly complex data types, including those generated by modern AI systems. Others are driven by the scale of contemporary datasets, which must often be stored, compressed, partitioned, or communicated under stringent resource constraints without sacrificing statistical efficiency in downstream inference. In such settings, the classical goal of uncovering low-dimensional structure becomes substantially more challenging. These challenges are further exacerbated when standard modeling assumptions fail, data are heterogeneous, or relevant structural representations must themselves be learned from the data.

The purpose of this paper is to provide a selective review of active research directions in high-dimensional statistics and articulate key challenges that are currently shaping the field. Given the rapid pace of methodological and theoretical developments, navigating the literature and identifying entry points into specific subareas has become increasingly difficult, particularly for researchers new to the field. This paper aims to synthesize recent advances across several representative topics, highlighting common themes, conceptual frameworks, and open problems, while pointing to foundational and contemporary references that serve as gateways to deeper engagement with particular research areas. By doing so, we hope to offer both a snapshot of the current landscape of high-dimensional statistics and a guide to several exciting directions. We focus our discussion on the following subareas of high-dimensional statistics:

\begin{itemize}
\item \textbf{Computational-statistical trade-offs:}
%\emph{When can underlying structure be uncovered in an algorithmically efficient manner?}
Computational-statistical trade-offs ask when the low-complexity structure that makes high-dimensional estimation possible (e.g., sparsity, low-rankness, clustering) can also be recovered by algorithmically efficient procedures. In many canonical problems, a sharp separation occurs between the information-theoretic threshold for detection/recovery and the signal strength or sample size required by known polynomial-time methods---phenomena that appear in sparse PCA, community detection, planted subgraph problems, tensor decomposition, and high-dimensional robust estimation. We discuss a growing toolkit for formalizing such gaps, including average-case reductions, statistical query lower bounds, sum-of-squares and low-degree predictions, and optimization landscape barriers, and highlight the complementary ways these frameworks illuminate what efficient algorithms can and cannot do. An important open question is how to make these results useful in practice, i.e., clarifying how additional structure, more data, or slightly higher computational budgets change what can be achieved; whether methods such as approximate message passing (AMP) are optimal among efficient algorithms or whether current lower bounds are simply not yet tight enough; and how these trade-offs appear in modern learning systems driven by optimization dynamics.

\item \textbf{Data integration:}
In many modern applications, including multi-omics studies and electronic health records in biomedicine, recommender systems, sensor networks, social and economic data analysis, computer vision, and natural language processing, data are collected from multiple sources or modalities. 
Different data sources may vary in scale, noise level, distribution, and overall data quality, and may share only partially overlapping signals, making it difficult to distinguish common structure from source-specific effects. Practical constraints further complicate data integration, including missing or partially matched observations, privacy and data-sharing restrictions, and limitations in communication, storage, or computation. Although data integration is inherently challenging, these difficulties are amplified in high-dimensional settings, where joint structures are harder to identify, and the scarcity of data requires more data-efficient learning methods. To address these issues, we highlight several active research directions: (i) multi-view learning, (ii) data integration under limited data access, (iii) integration of heterogeneous data, and (iv) distributed learning.

\item \textbf{High-dimensional asymptotics:} 
%\emph{How can asymptotic analysis help in the principled analysis of high-dimensional data?}
The performance of canonical frequentist statistical methods (e.g., empirical risk minimization, spectral methods, and maximum likelihood) on high-dimensional data has been investigated extensively over the past two decades. The high-dimensional asymptotic regime, where the number of samples and data points are both large, has emerged as a key framework in this endeavor. We highlight some key directions for future inquiry: (i) Bayesian inference in high-dimensional settings, (ii) optimal inference of low-dimensional functionals of high-dimensional parameters, and (iii) going beyond Gaussian designs under proportional asymptotics. Additionally, we note that many distinct  statistical methods are optimal under the traditional the minimax framework. This motivates the following question: Is it possible to compare competing methods by refining the minimax framework?  

\item \textbf{High-dimensional statistics and AI:} Artificial intelligence is transforming many areas of science and society. Modern AI systems rely on learning extremely large numbers of parameters from complex data, thereby challenging the assumptions of classical statistical theory and making high-dimensional statistics essential for understanding and improving AI models. Given the rapid pace and breadth of AI research, we focus on a few representative subtopics, such as fine-tuning, in-context learning, machine unlearning, and explainable AI, and discuss how tools from high-dimensional statistics can contribute to cross-disciplinary work.
\end{itemize}

While this paper covers a range of active research directions and highlights several important open problems in high-dimensional statistics, it does not aim to be an exhaustive survey of the field. The breadth and rapid evolution of high-dimensional statistics has given rise to many additional areas of active investigation that are beyond the scope of a single article. Notable topics such as uncertainty quantification, false discovery rate control, tensor analysis, and causal inference are only discussed briefly in Section \ref{ssec:discussion}; 
these omissions reflect choices of focus rather than a lack of importance.
%and we view them as complementary directions that will continue to play a central role in shaping the field. 

%\documentclass{article}
%\usepackage{fullpage}
%\usepackage{hyperref}
%\usepackage{amsfonts}
%\usepackage{graphicx}
%\usepackage{amsmath}
%\usepackage{amsthm}
%\usepackage{amssymb}
%\usepackage[english]{babel}
%\usepackage[letterpaper,top=2cm,bottom=2cm,left=3cm,right=3cm,marginparwidth=1.75cm]{geometry}
%\usepackage{xcolor}
%\usepackage{tikz}

%\newtheorem{example}{Example}

%\begin{document}

\section{Computational-Statistical Trade-offs}
\label{sec:trade-off}

%\begin{sloppypar}
\subsection{Introduction}

Computational-statistical trade-offs have emerged as a fundamental theme in modern high-dimensional statistics and machine learning. Classical statistical theory primarily concerns itself with information-theoretic limits and sample complexity; in contrast, a computational perspective recognizes that practical algorithms often face additional constraints beyond sample size. These constraints can manifest as gaps between what is statistically possible and what is computationally feasible in polynomial time.

The study of computational-statistical gaps has become increasingly important as researchers work with larger, more complex datasets. In many high-dimensional problems, regimes exist where efficient estimation is information-theoretically possible, yet no known polynomial-time algorithm can achieve the optimal statistical rate. Understanding these gaps helps practitioners choose appropriate algorithms and formulate problems effectively; for algorithm selection, if a problem exhibits a significant gap, practitioners may need to consider whether more sophisticated algorithms beyond standard methods are necessary, whether the problem parameters place them in a regime where efficient algorithms can succeed, or whether approximate solutions with weaker guarantees may be acceptable given computational constraints.
An important practical takeaway from lower bounds is understanding what will \emph{not} work: If lower bounds show that particular classes of algorithms cannot succeed, practitioners can focus their efforts elsewhere. For example, if spectral methods provably fail for a problem in certain regimes, practitioners should explore non-spectral alternatives, rather than spending time optimizing spectral methods.

Several canonical problems illustrate this phenomenon: In sparse principal component analysis, information-theoretic considerations allow detection and recovery despite comparatively weak signals, but a broad body of evidence suggests that polynomial-time methods require substantially stronger signals. Similar separations appear in network community detection, particularly under the stochastic block model, and are epitomized by the planted clique problem, where efficient algorithms are believed to require clique sizes far above the information-theoretic threshold. Computational-statistical gaps also arise in tensor decomposition and related higher-order inference tasks. More recently, such gaps have been identified and studied in various tasks in high-dimensional robust estimation \cite{diakonikolas2023algorithmic}.

Beyond problem-level separations that are agnostic to algorithm choice, it is also important to understand the limitations of specific algorithmic paradigms used by practitioners. Gradient-based methods, for instance, can fail in landscapes with unfavorable geometry, becoming trapped near poor local minima or lingering around saddle regions, even when statistically optimal recovery is possible. At the same time, modern observations about overparameterization and stochasticity suggest that optimization noise and implicit regularization can sometimes improve performance, highlighting a nuanced interplay between computation and statistics.
Spectral methods provide another widely used class of tools. Eigenvalue and singular value decompositions are natural in problems such as community detection and sparse principal component analysis. In many settings, spectral procedures achieve optimal or near-optimal performance in easier regimes, yet break down in harder regimes where more sophisticated approaches succeed. Characterizing where these transitions occur helps delineate when standard tools suffice and when one must move beyond them, motivating the need for general frameworks that formalize computational lower bounds and guide the development of new algorithmic methods.

\subsubsection{Frameworks for Computational-Statistical Gaps}

A variety of frameworks have been developed to study computational lower bounds in high-dimensional statistics, each capturing broad classes of algorithms and offering different forms of evidence for inherent computational barriers in the original problem. These approaches illustrate when statistical tasks become computationally intractable. We begin by reviewing the main frameworks and summarizing their respective strengths and limitations.

\paragraph{Reductions:}
Reductions from problems believed to be computationally hard provide a powerful method for establishing average-case computational hardness in statistical inference. The central idea is to show that if one could efficiently solve a given statistical problem, one could also efficiently solve another problem that is widely conjectured to be hard (on average). Such a reduction, once successfully established, demonstrates that any breakthrough algorithm for the target statistical task would necessarily imply a breakthrough for the source hard problem, thereby providing strong evidence of computational barriers. This technique has been used extensively to prove computational lower bounds in high-dimensional statistical problems such as sparse PCA, tensor PCA, tensor completion, submatrix detection, and community detection. The canonical conjectured hard problems involved include planted clique detection, hypergraphic planted clique, and random $k$-SAT, among others. Typically, the reduction maps instances of the hard problem to instances of the statistical estimation or testing task in a way that preserves the relevant structural properties, while ensuring that polynomial-time algorithms for the latter would yield polynomial-time algorithms for the former. Because these hardness assumptions are based on long-time failed algorithmic attempts and deep connections to worst-case complexity, reductions offer a flexible approach to deriving computational lower bounds.

%\subsection{Strengths and Limitations (Rishabh)}

\paragraph{Statistical Query Framework:}
In the statistical query (SQ) model \cite{kearns1998efficient}, the algorithm solving the inference task does not have unrestricted access to the dataset. Instead, it can only access the dataset by querying an oracle for the empirical average of a chosen scalar query function over the dataset.  The computational complexity of the algorithm is measured using the number of queries it makes. The responses to these queries are adversarially corrupted with noise of magnitude comparable to the random fluctuations in the empirical average expected under the model. The inference algorithm must be robust to these adversarial corruptions. The statistical query model exploits this robustness requirement to derive computational lower bounds. 

One of the key strengths of this framework is that the lower bounds obtained are very strong: For many inference problems with a computational-statistical gap, it has been shown that all SQ algorithms must make a super-polynomial or exponential number of queries to solve these problems in their hard phase \cite{feldman2017statistical,wang2015sharp,diakonikolas2017statistical}. Another attractive feature of the SQ model is that it uses robustness, a property likely necessary for any algorithm to be practically useful, to derive lower bounds. 

However, the robustness requirement in the SQ framework can be overly stringent. As a result, sharp SQ lower bounds can sometimes be too pessimistic, leading to incorrect predictions of the computational–statistical gap \cite{feldman2015complexity,li2019mean,dudeja2021statistical}. For example, in the tensor PCA problem of estimating an order-4, rank-1 signal tensor in $p$ dimensions from $n$ noisy samples, polynomial-query SQ algorithms require $n \gtrsim p^3$ samples for symmetric tensors and $n \gtrsim p^4$ samples for non-symmetric tensors.  In contrast, efficient spectral estimators succeed with $n \asymp p^2$ samples, implying that such estimators cannot be implemented within the SQ model. There have been some attempts to address this limitation \cite{feldman2015complexity}, but these fixes appear to be designed post hoc to ensure that the resulting SQ lower bounds capture a widely believed computational-statistical gap in a specific inference problem. 

%It would be interesting to further explore principled relaxations of the robustness requirements imposed in the SQ framework, which can correctly predict the computational-statistical gaps in a wide range of inference problems.

\paragraph{Landscape Analysis:}
A different framework provides evidence by studying popular optimization-based estimators for the problem, which estimate the unknown parameter by minimizing a cost function such as the negative log-likelihood. In the Overlap Gap Property (OGP) framework, one studies the sub-level sets of this cost function. For problems such as planted clique \cite{gamarnik2024landscape}, planted principal submatrix recovery \cite{gamarnik2021overlap,arous2023free}, and sparse linear regression \cite{gamarnik2022sparse}, it is known that the sub-level sets of the negative log-likelihood break into multiple disconnected clusters, and a randomly initialized algorithm typically starts at a  cluster which does not include the true signal. To move from one cluster to the other, an optimization algorithm must cross a region of high cost, which gradient descent and other local optimization methods avoid. This geometric property of the landscape is called an \emph{overlap gap property}, and rules out the success of popular greedy search methods such as gradient descent and certain Markov chains \cite{gamarnik2024landscape,gamarnik2021overlap,arous2023free}. 

One limitation of the OGP is that it is tailored to a particular cost function, and the class of algorithms it rules out may be limited. Even if a particular cost function exhibits an OGP, it is  possible that a different choice of cost function would be easy to optimize and result in a good estimator. Indeed, it is known that the OGP persists for certain natural cost functions (such as the negative log-likelihood) even in the ``easy phase'' of the planted clique problem, but a suitable modification of the cost function is easy to optimize and leads to a good estimator \cite{gamarnik2024landscape}.

\paragraph{Sum-of-Squares Hierarchy:}
The sum-of-squares (SoS) hierarchy provides a systematic framework for analyzing computational–statistical trade-offs in high-dimensional inference \cite{Parrilo2000,Lasserre2001}. Many statistical estimation and detection problems can be formulated as optimizing a polynomial subject to algebraic or combinatorial constraints. The SoS hierarchy constructs semidefinite programming (SDP) relaxations for such polynomial optimization problems. These SDP relaxations are parametrized by a degree parameter $d$, which simultaneously controls the strength of the relaxation and the algorithmic cost: a degree-$d$ SoS relaxation can typically be solved in time $n^{O(d)}$. Consequently, low-degree SoS relaxations correspond to polynomial-time or mildly subexponential algorithms, while higher degree $d$ interpolates toward brute-force search. This makes SoS a proxy for ``all efficient algorithms," in the sense that many classical algorithms, such as spectral methods, SDPs, local search procedures, and certain message-passing algorithms, can often be reformulated as or dominated by low-degree SoS algorithms \cite{barak2014sum}.

To establish computational hardness using the SoS hierarchy, one typically shows that low-degree SoS fails to solve a given statistical problem in signal-to-noise regimes that are conjectured to be computationally difficult. Some examples in the literature include sparse PCA \cite{MekaPotechinWigderson2015}, tensor PCA \cite{HopkinsSchrammShiSteurer2015}, submatrix detection and planted sparse structures \cite{hopkins2018integrality}, and refuting random constraint satisfaction problems. While SoS computational lower bounds are often technically intricate to prove, the SoS hierarchy is regarded as one of the strongest known algorithmic frameworks, covering a wide range of existing algorithmic approaches. Therefore, showing the failure of low-degree SoS relaxations on high-dimensional statistical problems is viewed as compelling evidence of an inherent computational barrier.

\paragraph{Low-Degree Framework:}
The low-degree polynomial framework is also seen as a ``lite version" of the sum-of-squares hierarchy and is based on a conjecture of Hopkins \cite{hopkins2018statistical}
%and inspired by the lower bounds for the sum-of-squares hierarchy \cite{barak2019nearly,hopkins2017bayesian,hopkins2017power}, a powerful family of semi-definite programs studied in theoretical computer science. 
At a high level, the conjecture states that, for many natural statistical testing and estimation problems involving a $p$-dimensional unknown parameter, the class of estimators that are degree-$D = O(\log(p))$ polynomials in the observed dataset is as powerful as all computationally efficient (polynomial-time) estimators. Thus, the low-degree framework provides evidence of a computational-statistical gap by proving that all degree-$D = O(\log(p))$ (or even higher degree) polynomials fail to solve the testing or estimation problem in the ``hard'' phase. 

There are two main variants of the low-degree method: one for testing \cite{hopkins2018statistical,kunisky2019lowdegree} and one for estimation \cite{schramm2020computational}. In testing problems, the framework shows that all low-degree tests fail to separate the null and alternative hypotheses, i.e., the absolute difference between the expected value of the test under the null and the alternative is negligible relative to the test's standard deviation.  For estimation problems, the low-degree framework aims to obtain lower bounds on the minimum mean-squared error achievable by degree-$D$ polynomial estimators. This lower bound shows that these estimators are no better than a trivial constant estimator (which does not use the observed data) in the hard phase. The low-degree framework has recently been extended to obtain lower bounds for algorithms that are not polynomials of the observed dataset. Examples include tests obtained thresholding a low-degree polynomial (known as polynomial threshold functions) \cite{diakonikolas2025ptf}, or applying an arbitrary transformation to a low-degree polynomial \cite{hsieh2026rigorous},  and low-coordinate-degree functions---represented as linear combinations of functions depending on only a few coordinates of the dataset \cite{kunisky2025low}.

A key strength of the low-degree framework is its versatility and broad applicability. It has been applied to various problems and consistently identifies the correct computational-statistical gap \cite{luo2024computational,lyu2023optimal,luo2020tensor,luo2024tensor,sohn2025sharp,jinphase,ding2019subexponential,d2020sparse,loffler2020computationally,bandeira2020computational,arpino2023statistical}. However, the framework uses the degree of a polynomial as a proxy for computational cost. Sometimes, high-degree polynomials have additional structure that permits efficient computation. Some works have exploited this limitation to construct efficient algorithms for statistical problems that the low-degree framework predicts to be computationally hard \cite{zadik2022lattice,diakonikolas2022non,buhai2025quasi,koehler2022reconstruction}. \\

\paragraph{Summary:} Each of the frameworks described above has distinct characteristics. The statistical query framework is quite general and captures many natural algorithms, but the lower bounds can sometimes be circumvented by algorithms that use the data in non-SQ ways. The low-degree framework provides refined predictions about algorithmic performance but requires careful problem-specific analysis. Sum-of-squares lower bounds are algorithmically meaningful, but proving such bounds can be technically challenging. Reductions provide strong evidence of hardness when they exist, but finding appropriate reductions can be difficult and may depend on unproven complexity assumptions. An important observation is that different frameworks can give different predictions for the same problem. Understanding when frameworks agree or disagree, and what these differences reveal about the problem structure, remains an active area of research.

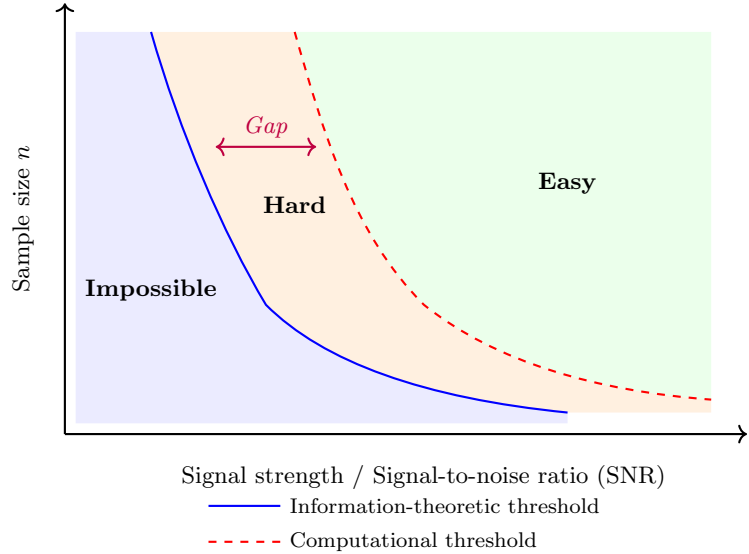
\begin{figure}[ht]
\centering
\begin{tikzpicture}[scale=0.95]
  % --- Shaded regions (behind everything) ---
  % Impossible region: left of blue curve
  \fill[blue!8] (0.15,0.15) -- (0.15,5.6)
    -- (1.2,5.6) .. controls (1.6,4.2) and (2.2,2.8) .. (2.8,1.8)
    .. controls (3.6,1.0) and (5.0,0.5) .. (7.0,0.3)
    -- (7.0,0.15) -- cycle;
  % Hard region: between the two curves
  \fill[orange!12] (1.2,5.6)
    .. controls (1.6,4.2) and (2.2,2.8) .. (2.8,1.8)
    .. controls (3.6,1.0) and (5.0,0.5) .. (7.0,0.3)
    -- (9.0,0.3)
    -- (9.0,0.48)
    .. controls (7.4,0.6) and (6.0,1.0) .. (5.0,1.8)
    .. controls (4.0,2.8) and (3.6,4.2) .. (3.2,5.6)
    -- cycle;
  % Easy region: right of red curve
  \fill[green!8] (3.2,5.6)
    .. controls (3.6,4.2) and (4.0,2.8) .. (5.0,1.8)
    .. controls (6.0,1.0) and (7.4,0.6) .. (9.0,0.48)
    -- (9.0,5.6) -- cycle;

  % --- Axes ---
  \draw[->, thick] (0,0) -- (9.5,0);
  \draw[->, thick] (0,0) -- (0,6.0);

  % --- Axis labels (properly placed) ---
  \node[below] at (5.0,-0.3) {\small Signal strength / Signal-to-noise ratio (SNR)};
  \node[rotate=90, above] at (-0.3,3.0) {\small Sample size $n$};

  % --- Information-theoretic threshold curve (solid blue) ---
  \draw[thick, blue] (1.2,5.6) .. controls (1.6,4.2) and (2.2,2.8) .. (2.8,1.8)
    .. controls (3.6,1.0) and (5.0,0.5) .. (7.0,0.3);

  % --- Computational threshold curve (dashed red) ---
  \draw[thick, red, dashed] (3.2,5.6) .. controls (3.6,4.2) and (4.0,2.8) .. (5.0,1.8)
    .. controls (6.0,1.0) and (7.4,0.6) .. (9.0,0.48);

  % --- Region labels (well inside each region) ---
  \node at (1.2,2.0) {\small \textbf{Impossible}};
  \node at (3.2,3.2) {\small \textbf{Hard}};
  \node at (7.0,3.5) {\small \textbf{Easy}};

  % --- Gap arrow (between curves, clear of both) ---
  \draw[<->, thick, purple] (2.1,4.0) -- (3.5,4.0);
  \node[above, purple] at (2.8,4.0) {\small \textit{Gap}};

  % --- Legend ---
  \draw[thick, blue] (2.0,-1.0) -- (3.0,-1.0);
  \node[right] at (3.0,-1.0) {\footnotesize Information-theoretic threshold};
  \draw[thick, red, dashed] (2.0,-1.5) -- (3.0,-1.5);
  \node[right] at (3.0,-1.5) {\footnotesize Computational threshold};
\end{tikzpicture}
\caption{Schematic phase diagram illustrating the computational-statistical gap. The solid blue curve marks the information-theoretic threshold below which no algorithm can succeed. The dashed red curve marks the computational threshold below which no known polynomial-time algorithm succeeds. The shaded region between the two curves is the \emph{hard} regime.}
\label{fig:phase-diagram}
\end{figure}

Figure~\ref{fig:phase-diagram} illustrates the general structure of a computational-statistical gap using a schematic phase diagram. In many canonical problems, the parameter space is partitioned into three regimes: an \emph{impossible} regime where no algorithm can succeed regardless of computational resources, a \emph{hard} regime where information-theoretically optimal estimators exist but no known polynomial-time algorithm succeeds, and an \emph{easy} regime where efficient algorithms achieve the statistical goal.

\subsubsection{Overview}

The remainder of this section develops three themes in detail: (i) Section \ref{ssec:hierarchyhardness} discusses the hierarchy of hardness results and connections between different frameworks; (ii) Section \ref{ssec:realisticmodel} addresses bridging theory and practice; and (iii) Section \ref{ssec:blessing:barrier} explores the counterintuitive phenomenon whereby computational barriers can yield inferential blessings.

\subsection{Hierarchy of hardness results}\label{ssec:hierarchyhardness}

\paragraph{Background:} Most existing computational lower bounds in statistics involve reducing well-studied problems in theoretical computer science, which are generally accepted to be computationally hard, to the statistical problem being studied. However, types of ``hard" problems may still differ in terms of difficulty, i.e., in terms of the likelihood that the conjectured hardness holds, or in terms of their feasibility if one is willing to use a slightly larger computational budget. Computational bounds should always be approached with a critical eye. Hardness results are traditionally derived in a worst-case sense, e.g., in the very definition of NP-hardness. For practical purposes, however, one may be more interested in understanding the runtime of algorithms applied to ``typical" inputs, drawn from a certain distribution that might place negligibly low probability on such worst-case inputs.

\paragraph{Recent work:} The intriguing book~\cite{roughgarden2021beyond}, written by various contributors from the theoretical computer science community, explains the benefits of taking a more nuanced, average-case approach to a variety of classical and contemporary problems, including sparse recovery, stochastic block models, robust estimation, and matrix factorization. In particular, the authors emphasize that the theoretical analysis of algorithms is beneficial from the point of view of (1) performance prediction, (2) identification of optimal algorithms, and (3) development of new algorithms, for which a paradigm shift from studying worst-case to  average-case performance of algorithms could be  fruitful.

To make these concepts  more concrete, consider the classical problem of high-dimensional linear regression. Under well-established conditions such as a restricted eigenvalue condition, one can derive both good statistical accuracy of the global optimum of a regularized least squares objective and fast convergence of convex optimization algorithms~\cite{agarwal2012fast}. When it comes to verifying these conditions, however, it is widely accepted in statistical theory that the restricted eigenvalue condition holds \emph{with high probability} over a certain data distribution. Thus, there is a clear benefit in statistical analysis to consider \emph{average-case} performance of algorithms. Relatedly, \cite{paquette2021sgd} develop an average-case analysis of stochastic gradient descent, further illustrating how distributional modeling can lead to more refined predictions of algorithmic behavior beyond worst-case guarantees.

\paragraph{Future directions:}
The work \cite{berthet2013complexity} derives a computational lower bound for sparse principal component detection in terms of the hardness of a hypothesis test for detecting a planted clique in an Erd\H{o}s-Renyi random graph; the latter problem is widely conjectured to be hard when the size of the planted clique is sub-quadratically small. Planted clique detection in Erd\H{o}s-Renyi graphs may be classified as an ``average-case complexity" result. Beyond the planted clique problem, higher-order analogs motivate the hypergraphic planted clique (HPC) as a unifying hardness hypothesis; \cite{luo2020open} formally posed the HPC open problem and called for equivalence results between HPC and the standard planted clique. Exploring this connection further is a promising direction for future work.

In contrast, \cite{zhang2014lower} establish computational lower bounds for sparse linear regression by comparison to a constraint satisfaction problem called the exact 3-set cover problem, which is known to be NP-hard. Since NP-hardness is a standard notion in worst-case complexity theory, the validity of the computational lower bounds in \cite{zhang2014lower} thus merely depend on the validity of the conjecture \textbf{NP $\not\subset$ P/poly}, where the latter set corresponds to the class of problems solvable in polynomial time by a Turing machine with an advice string of polynomial length (hence a superset of $P$). Further exploring such nuances between hardness assumptions required to derive statistical conclusions could be another fertile area of research.

Another line of recent work \cite{regev2009lattices, bruna2021continuous} analyzes problems in learning theory, such as density estimation for Gaussian mixtures, which are shown to be computationally hard under conjectures regarding the difficulty of solving certain learning problems under \emph{quantum} algorithms. The shift of interest to hardness under quantum computation rather than classical computation seems very relevant for the future. We note that a better understanding of the sorts of lattice-based algorithms analyzed in this line of work, and their hardness under different paradigms of computation and potential applications, could be quite beneficial to statisticians. Indeed, while these problems are actively studied in cryptography, their ramifications for statistical problems are only beginning to be realized \cite{bun2023stability, gupte2022continuous}.

Understanding the relationships between different frameworks for computational-statistical gaps remains an important open problem. The frameworks reviewed above---reductions, statistical queries, landscape analysis, sum-of-squares, and low-degree polynomials---have largely been developed and applied independently, yet they often target the same underlying phenomena. Some progress has been made: for instance, in certain settings, low-degree polynomial lower bounds imply statistical query lower bounds, and vice-versa \cite{brennan2021statistical}. Similarly, connections between the sum-of-squares hierarchy and the low-degree framework have been explored, as the low-degree conjecture was originally inspired by SoS lower bounds \cite{hopkins2018statistical,barak2019nearly}. Connections between heuristics from statistical physics and low-degree and SQ lower bounds have also been explored \cite{bandeira2022franz,chen2025optimized,tsirkas2026monotonicity}. However, the precise relationships between these frameworks are not fully understood. Can we identify natural statistical problems where one framework predicts hardness but another does not? How do algorithm-specific lower bounds (e.g., from landscape analysis) relate to algorithm-agnostic frameworks (e.g., low-degree or SoS)? Developing unified frameworks that capture the strengths of multiple methods remains an important goal.

%%%%%

\subsection{Bridging theory and practice}
\label{ssec:realisticmodel}

\paragraph{Background:}
Most theoretical work on computational--statistical gaps has focused on relatively simple data-generating processes. These models often assume independent, homogeneous observations without outliers or contamination, whereas such assumptions are frequently violated in practice.  Another limitation is that much of the literature considers simple parametric models, while real-world applications often rely on far more complex models, such as neural networks.  

%\paragraph{\textcolor{blue}{Recent work:}}

\paragraph{Future directions:} Understanding computational-statistical gaps in more realistic models, such as settings with dependent data or structured designs, remains an important open area. Questions about time series, spatial data, network-structured observations, heterogeneous data, outliers, contamination, missing data, or measurement error are largely unexplored. Tools from random matrix theory, statistical physics methods, and concentration inequalities adapted to high-dimensional settings play important roles in analyzing these problems. Rigorizing predictions from statistical physics remains an important direction, with connections to the low-degree polynomial framework and other mathematical methods.

Strengthening connections to modern machine learning practice is also an important goal. How do computational-statistical trade-offs manifest in deep learning? What can the theory reveal about phenomena such as double descent and benign overfitting? Can we use ideas from computational-statistical gaps to design better architectures and training procedures?
In the context of neural networks, understanding computational barriers can inform architecture choices. Research has shown that for certain problems such as single-index and multi-index models, the depth and width of neural networks can affect whether computational barriers are overcome \cite{abbe2023sgd,damian2024computational,bietti2022learning,damian2022neural}.  Some problems may benefit from additional layers or specific architectural features that help avoid computational barriers.

Finally, a major limitation of existing work on computational-statistical trade-offs is the reduction to algorithms that are theoretically polynomial-time, but not practically implementable. For instance, while sum-of-squares methods have gained prominence for yielding polynomial-time algorithms with (nearly) optimal statistical accuracy, the space complexity of a degree-$d$ sum-of-squares optimization problem with $n$ variables scales as $O(n^d)$. In statistical applications, one typically encounters large $n$ and moderate $d$; as noted in \cite{barak2014sum}, ``even $d = 6$ and $n = 100$ would take something like a petabyte of memory \dots\ [we] had a hard time executing a program with $n = 16$ and $d = 4$ on a Cornell cluster'' \cite{barak2014sum}. It is an active area of research in optimization to determine ways to speed up solvers for sum-of-squares programs, recognizing their widespread adoption as a theoretical tool in other areas of computational mathematics and statistics~\cite{ahmadi2017improving}. Thus, for many of the aforementioned problems, there is still a dearth of algorithms in the literature which achieve optimal statistical performance, while being truly practically implementable.

%%%%%

\subsection{Blessing of computational barriers}\label{ssec:blessing:barrier}

\paragraph{Background:}

A recurring theme in recent work is a counterintuitive ``blessing of computational barriers'': regimes in which the additional signal strength or sample size required for tractable computation incidentally yields simpler or sharper inference.

\paragraph{Recent work:}

For low-rank tensor models, for example, feasible algorithms typically need stronger conditions than the information-theoretic limit; those very conditions can deliver asymptotic normality and valid uncertainty quantification without debiasing, a phenomenon formalized for tensor PCA/regression in \cite{xia2022inference}. Likewise, in tensor-on-tensor regression, over-parametrization interacts favorably with Riemannian optimization: moderate rank over-parametrization can be essentially ``cost-free'' in sample size while aligning with computationally feasible procedures, and low-degree arguments reveal genuine statistical–computational gaps in related scalar-on-tensor problems \cite{luo2024tensor}.

\paragraph{Future directions:}

An interesting direction is to systematize when and how such barriers become blessings: characterize signal-to-noise/sample thresholds that simultaneously (i) ensure polynomial-time estimation and (ii) yield valid inferential procedures across a broader class of high-order models. For instance, one could ask whether the blessing phenomenon extends to broader families of structured estimation problems, such as higher-order tensor decompositions with heterogeneous noise, or to settings involving dependent or incomplete observations. A related but distinct perspective is the ``blessing of high-dimensionality'' \cite{li2018embracing}, which observes that increasing the ambient dimension can sometimes simplify rather than complicate inference---for example, by improving the conditioning of certain optimization landscapes or enabling more effective regularization. Understanding how this phenomenon interacts with computational barriers, and whether the additional structure imposed by high dimensionality can simultaneously ease both computational and inferential challenges, is a promising avenue for future investigation.

\subsection{Further open directions}

We close with a brief summary of additional open directions that cut across the themes of this section:
\begin{itemize}
\item \textbf{Unifying frameworks.} Understanding when different frameworks (statistical queries, low-degree polynomials, sum-of-squares, reductions) agree or disagree in their predictions, and developing unified theories that capture the strengths of multiple approaches, remains an important open problem.
%\item \textbf{Realistic data models.} Extending computational-statistical gap analyses to settings with dependent observations, outliers, heterogeneous data, or missing data is largely unexplored and represents a major opportunity for future work.
\item \textbf{Gaps at the level of constants.} Most existing work characterizes computational-statistical gaps in terms of rates. However, for problems such as estimating a rank-1 spike in a Gaussian matrix, the gap between information-theoretically optimal and computationally efficient methods manifests at the level of sharp constants---specifically, the precise SNR threshold. It is widely believed that approximate message passing (AMP) algorithms \cite{donoho2009message, feng2022unifying} achieve the optimal SNR threshold among efficient methods, but rigorous evidence remains limited: current low-degree lower bounds only rule out constant-degree polynomials \cite{montanari2025equivalence} rather than the logarithmic-degree polynomials that would provide stronger evidence. Understanding when AMP is optimal among efficient procedures, and extending constant-level gap analyses more broadly, is an important open direction. See \cite{sohn2025sharp} for recent progress in this direction.
%\item \textbf{Connections to modern machine learning.} Formalizing how computational-statistical trade-offs manifest in deep learning---including phenomena such as double descent and benign overfitting---and using these insights to guide architecture design and training procedures is an important open direction.
\item \textbf{Sample size, computation, and privacy.} Understanding whether additional samples can compensate for computational constraints, and how privacy requirements interact with both statistical and computational limitations, remains an emerging and largely open area.
\end{itemize}

%%%%%

%\subsection{conclusion \textcolor{red}{arian: maybe we can remove this section?}}

%Computational-statistical trade-offs represent a fundamental aspect of modern high-dimensional statistics. Understanding these trade-offs requires synthesizing ideas from statistics, computer science, optimization, and information theory. Significant progress has been made in developing frameworks for characterizing computational barriers and identifying specific problems where such barriers arise. Nevertheless, many questions remain open. Future work bridging theory and practice, connecting different frameworks, and extending results to more realistic settings will be essential for making this theory maximally useful for practitioners.

%\bibliographystyle{plain}
%\bibliography{references-trade-off}
%\end{sloppypar}

%\end{document}

%\documentclass{article}
%\usepackage{fullpage}
%\usepackage{hyperref}
%\usepackage{lipsum}
%\usepackage{amsfonts}
%\usepackage{graphicx}
%\graphicspath{ {images/} }
%\usepackage{epstopdf}
%\usepackage{dsfont}
%\usepackage{algorithmic}
%\usepackage{mathtools}
%\usepackage{empheq}
%\usepackage{amsfonts}
%\usepackage{amsgen,amsthm,amsmath,amstext,amsbsy,amsopn,amssymb,subfigure,stmaryrd,mathabx,mathrsfs}
%\usepackage[dvips]{graphicx}
%\usepackage{comment}
%\usepackage[font=footnotesize,labelfont=bf]{caption}
%\usepackage{xcolor}

%\usepackage{blkarray}

%\usepackage{url}
%\usepackage{longtable}
%\usepackage{mathtools}
%\usepackage{multirow}
%\usepackage{array}
%\begin{document}

\section{Data integration}

\subsection{Introduction}

Data integration is a key challenge in high-dimensional statistics, especially when analyzing data from biomedical applications \cite{tseng2015integrating, margolis2014national, weinstein2013cancer}. Existing data integration approaches can be broadly divided into two categories, \emph{horizontal} and \emph{vertical} integration \cite{richardson2016statistical}. 
In \emph{horizontal data integration}, data are integrated across multiple studies/sources. Specifically, consider $m$ studies containing data with the same set of $p$ covariates: $\bm{X}^{(j)} \in \mathbb{R}^{n_j \times p}$, for $j = 1, \ldots, m$. 
A classical example of horizontal data integration is 
\emph{meta-analysis} \cite{begum2012comprehensive, van2002advanced, normand1999meta, evangelou2013meta}, where raw data or summary statistics are integrated across multiple studies to improve power and reproducibility  \cite{nalls2014large, shen2010meta}.

In \emph{vertical data integration}, the datasets contain different sets of variables collected on the same set of subjects: $\bm{X}^{(k)} \in \mathbb{R}^{n \times p_k}$, for $k = 1, \ldots, q$. In this case, data integration aims to obtain a more comprehensive perspective by  studying different factors associated with the outcome. An important application is the integration of multiple types of \emph{omics} data in association with a disease or phenotype \cite{gomez2014data, curtis2012genomic, mo2013pattern}, where \emph{multi-view analysis} methods \cite{li2018review} are commonly used to integrate across omics data types. 

In what follows, we first discuss multi-view learning, one of the most well-known vertical data integration settings, and devote the rest of the section to emerging problems in horizontal data integration. %\textcolor{red}{Everyone please categorize your section accordingly.}
% Sections \textcolor{red}{X,} , \ref{SecDistributed} focus on horizontal data integration, whereas Sections \textcolor{red}{Y} focus on vertical data integration and Sections \textcolor{red}{sec:limiteddata} discuss aspects of both.

% NOTE: For each of the sections, need to discuss (i) stat problems; (ii) algorithms; (iii) applications and (iv) \emph{future directions}

\subsection{Multi-view learning}
\label{SecMulti}

\paragraph{Background:}
Multi-view learning focuses on integrating multiple data modalities observed for the same set of subjects. For simplicity, consider two data modalities, $\bm{X}^{(1)} \in \mathbb{R}^{n \times p_1}$ and $\bm{X}^{(2)} \in \mathbb{R}^{n \times p_2}$, corresponding to measurements of dimensions $p_1$ and $p_2$. 
Multi-view learning can be performed in an \emph{unsupervised} or \emph{supervised} fashion. If the diﬀerent data views share underlying relationships, integrating data across modalities can offer advantages over separate analyses of individual data modalities.

Unsupervised multi-view learning commonly focuses on identifying joint patterns in multiple data modalities. Examples include clustering $n$ individuals based on the joint signal in $\bm{X}^{(1)}$ and $\bm{X}^{(2)}$ \cite{bickel2004multi, chao2021survey, yang2018multi}, or learning joint (low-dimensional) representations \cite{li2018survey, wang2015deep}. 
The joint structure captures shared information in multi-view data. In contrast, the individual structure from each data view captures information that is unrelated to the joint structure. For two data views, a simple model capturing the joint and individual structures posits that \cite{lock2013joint}
\begin{align*}
    \bm{X}^{(1)} &= \bm{J}^{(1)} + \bm{A}^{(1)} + \bm{\varepsilon}^{(1)}, \\
    \bm{X}^{(2)} &= \bm{J}^{(2)} + \bm{A}^{(2)} + \bm{\varepsilon}^{(2)},   
\end{align*}
where $\bm{J} = [\bm{J}^{(1)}; \bm{J}^{(2)}]^\top$ represents the joint structure,  $\{\bm{A}^{(k)}\}_{k=1,2}$ captures individual structure, and $\mathbb{E}(\bm{\varepsilon}^{(k)}) = 0$.  
Various machine learning approaches for unsupervised multi-view learning try to learn these joint and individual structures and/or use them to perform downstream tasks such as clustering or dimension reduction. These include multi-view factor analysis and dimension reduction \cite{lock2013joint, hou2010multiple, kan2015multi}, canonical correlation analysis (CCA) \cite{rupnik2010multi}, non-negative matrix factorization (NMF) \cite{liu2013multi, wang2017diverse}, kernel-based methods \cite{tzortzis2012kernel}, network-based approaches \cite{kan2016multi, wang2019gmc}, and variational auto-encoders (VAEs) \cite{zhao2024deep, zhang2019integrated}; see, e.g., \cite{li2018review} for a review. 

Supervised multi-view learning also focuses on integrating the signal across multiple views. Depending on the application, supervised multi-view learning methods may focus on  predicting the outcome, or identifying  subsets of features in the the multi-view data that are associated with the outcome. 

\paragraph{Recent work:}
Multi-view learning methods, both supervised and unsupervised, can be broadly categorized based on the stage at which data fusion occurs: early, late, or intermediate~\cite{lanckriet2004statistical}. 

Early fusion methods form a single joint representation of the multi-view data and use this representation in downstream supervised or unsupervised learning \cite{hou2020multiview}. 
The joint representation can be simply constructed by concatenating the features of the multi-view datasets, e.g., by forming a concatenated matrix $[\bm{X}^{(1)}|\bm{X}^{(2}] \in \mathbb{R}^{p_1+p_2}$. Alternatively, early fusion can be performed by first inferring a (lower-dimensional) representation for each data modality---e.g., using PCA or VAE---and then performing supervised or unsupervised learning by combining the representations. 
% In the first approach, the concatenated feature matrix $\mathcal{X}$ is used in regression/classification, 
Despite its simplicity, the first approach can lead to larger feature sets and exacerbated multi-collinearity among features.
Therefore, multi-view data integration using lower-dimensional representations has gained more popularity in high-dimensional settings. Several approaches have also been proposed for learning joint representations among multi-view data \cite{hou2010multiple, kan2015multi}, capturing both the common structure of multiple data views and structures unique to each data view. 
The common structure can be used for identifying the joint clustering of observations based on multi-view data, whereas the unique structures can be used to find associations between each data view and the outcome. 
Although most joint representation methods---including JIVE \cite{lock2013joint}, MOFA \cite{argelaguet2018multi}, and SLIDE \cite{gaynanova2019structural}---focus on unsupervised learning, supervised learning methods have also been proposed for predicting the outcome based on joint representations. Examples of approaches for supervised learning include inferring linear regulatory modules between data types \cite{zhu2016integrating}, extracting joint and individual structures relevant to the outcome \cite{palzer2022sjive, wang2024joint}, and modeling each data view's low-rank signal using factor regression \cite{li2022integrative}. 

Late fusion methods perform supervised or unsupervised learning on each data view and then integrate the results. This includes integrating cluster memberships or predictions learned from each data view by averaging or majority voting \cite{liu2018late, wang2019multi, morvant2014majority}. As expected, intermediate fusion methods fall somewhere between early and late fusion. A recent example is collaborative learning \cite{gross2015collaborative, ding2022cooperative},
which combines early and late fusion using a tuning parameter. 

\paragraph{Future directions:}
The computational and theoretical considerations discussed elsewhere in this article also arise in the context of multi-view data integration. For instance, 
joint dimension reduction methods give rise to similar asymptotic questions as those discussed for dimension reduction for a single data view in Section~\ref{sec:pcaexample}. 
Moreover, identifying joint and individual (lower-dimensional) representations also requires additional identifiability assumptions \cite{andreou2019inference}. Most supervised joint representation learning approaches (including those discussed above) have been developed in the context of linear regression and do not apply to generalized linear models or methods for time-to-event data. 

Several statistical issues unique to multi-view learning also need to be addressed. First, uncertainty quantification and valid statistical inference for multi-view learning has not received adequate attention. Second, common representation learning approaches for multi-view learning, especially those involving deep neural networks such as VAEs, often provide limited interpretability. Third, as discussed in Section~\ref{SecMissing}, multi-view data may involve significant  missingness. For instance, biomedical applications often involve subsets of individuals for whom one or more data views are not observed. On the one hand, ignoring available observations with missing data in these \emph{block missingness} scenarios can lead to significant loss of information/power. On the other hand, traditional imputation strategies may not be effective in such block missing settings. Leveraging generative AI for synthetic data (see Section~\ref{sec:limiteddata}) or developing methods that can integrate partially matched multi-view data (as in \cite{yu2020optimal}) is an important direction for future research. 

%%%%%
\subsection{Integration in limited data access settings}
\label{sec:limiteddata}

%\subsubsection{Learning from summary statistics}
%\label{SecSummary}

%\textcolor{red}{It would be good to better distinguish this discussion from that of Section~\ref{SecDistributed}. In Section~\ref{SecDistributed}, we focus on homogeneous data. Is it true that summary statistics focuses on heterogeneous settings? Or is there indeed some overlap between the two topics? \\
%Re: There is certainly some connection especially with respect to meta-analysis where homogenous data is assumed and heterogeneity evaluated, yet other tasks like variable selection do not make any assumptions and allow for shared and distinct signals.}

Many applications of high-dimensional data integration are characterized by limited data access due to privacy restrictions, proprietary concerns, or technical barriers, since data are often stored in secure computing environments or data silos that prohibit simple data aggregation. These scenarios require the creation of novel data structures and analytic strategies. In the following, we focus on statistical methodologies for summary-level and synthetic data.

\paragraph{Background and recent work:}
%\subsection{Background}
%What is summary-level data?
Summary-level data refers to summary statistics such as effect sizes, standard errors, and $p$-values, derived from individual-level data. %The use of summary-level data is motivated mainly, but not exclusively, by health data analytics, where data are confidential and sensitive and therefore cannot be easily shared. In addition to privacy and regulatory constraints, summary statistics are often used because individual studies may have limited sample sizes, making it difficult to reliably estimate more complex or high-dimensional quantities. In such settings, simpler, aggregated summaries can be estimated more accurately and then combined across multiple sources to improve statistical power and stability. 
Key examples include genetics and genome-wide association studies \cite{Pasaniuc2017} or electronic health records \cite{Albers2018}. For example, in genome-wide association studies, each genetic variant is regressed against the outcome of interest in a regression model accounting for known confounders and technical covariates. The resulting summary statistic for the association strength of a single genetic variant captures the regression coefficient, standard error, and significance. Apart from practical considerations of data sharing, working on summary-level data allows one to maximize sample size and power and offers flexibility to integrate data from different sources to draw conclusions about relationships between variables. This is infeasible on an individual level, since not all variables may have been measured in a single dataset.

\begin{figure}[t]
\begin{center}
\includegraphics[width=6in]{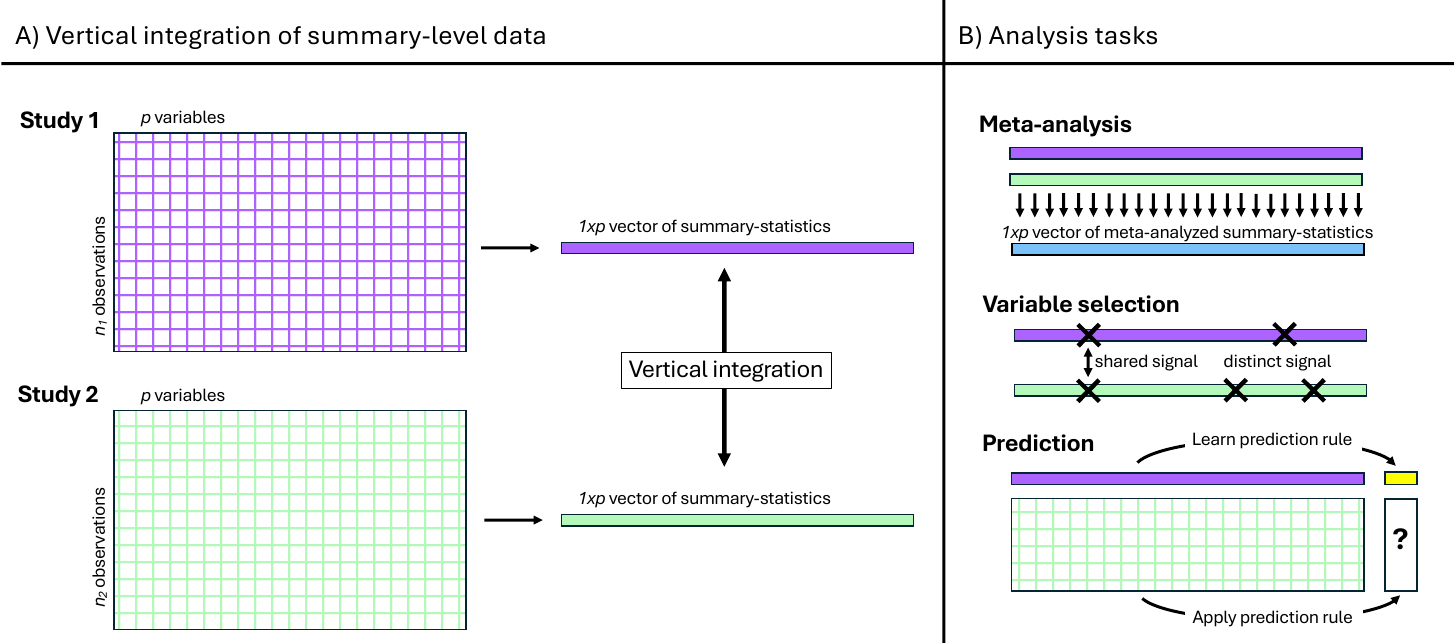}\label{fig:summary}
\caption{(a) Schematic illustration of data integration from summary level data. (b) Illustration of analysis task integrating summary-level data.}
\end{center}
\end{figure}

%\paragraph{Recent work}
Statistical tasks for data integration of summary-level data mostly involve horizontal integration, where data $\bm{X}^{(j)} \in \mathbb{R}^{n_j \times p}$, for $j = 1, \ldots, m$, are summarized or compressed into data vectors $\bm{s}^{(j)} \in \mathbb{R}^{1 \times p}$ of summary statistics, also referred to as \emph{proxy data} \cite{Li2022}. Information on the $p \times p$ covariance matrix $\bm{\Sigma}$ between the covariates is available in many applications. As shown in Figure~\ref{fig:summary}, there are three main tasks when working on summary-level data: Meta-analysis, variable selection, and prediction.
%\begin{itemize}
%\item \textit{Meta-analysis}: 
The aim of meta-analysis is to combine and aggregate information across $m$ studies and quantify heterogeneity between studies. Although meta-analysis was initially developed for a single variable of interest, it has been extended to handle a high-dimensional set of covariates while accounting for the correlation among covariates \cite{Liyanage2024}. Another aspect of horizontal data integration is the prospect of boosting power. This is particularly relevant when data from a single target population are scarce, leading to low power for the target population, but data from ($m-1$) other ancillary populations related to the target population are available \cite{Mason2025}. (See also Section~\ref{SecTransfer} for connections with transfer learning.)
% \textcolor{red}{This might connect to Section~\ref{SecTransfer}? Re: Indeed this links beautifully with the transfer learning. Thanks so much for the comment.}
% %\item \textit{Variable selection}: 
While univariate variable selection on summary-level data considers a $p$-dimensional vector of summary statistics $\bm{s} \in \mathbb{R}^{1 \times p}$ and aims to identify or select important features for one trait of interest \cite{Gao2020}, data integration of summary-level data allows for the detection of shared signals across multiple datasets. A common application in genetics is the detection of genetic variants associated with more than one trait, also known as pleiotropy \cite{Giambartolomei2014}. In this case, each of the $m$ traits to be analyzed is presented as vector of summary statistics $\bm{s}^{(j)} \in \mathbb{R}^{1 \times p}$, for $j = 1, \ldots, m$, measured on the same $p$ covariates.
%\item \textit{Prediction}: 
Prediction algorithms can also be trained on summary-level data alone for single- \cite{Zhang2021} or multi-task learning \cite{Knight2023}. Regularized regression \cite{Knight2023}, Bayesian approaches allowing for covariate-specific priors \cite{Zhang2021}, or ensemble learning \cite{Chen2024} have been used in this setting. 
%\end{itemize}

Synthetic data generation offers an alternative paradigm for both horizontal and vertical data integration. %when direct pooling of raw observations is infeasible due to privacy restrictions, proprietary concerns, or technical barriers to data sharing. 
In horizontal integration, consider $m$ studies with data $\bm{X}^{(j)} \in \mathbb{R}^{n_j \times p}$ drawn from distributions $\{\mathcal{P}_j\}_{j = 1, \dots, m}$. Rather than pooling the raw data, synthetic data methods generate $\bm{\widetilde{X}}^{(j)} \sim \widehat{\mathcal{P}}_j$ from each study, where $\widehat{\mathcal{P}}_j$ is an estimate of $\mathcal{P}_j$ based on the original data, and then combine these synthetic datasets to perform meta-analysis or joint estimation. In vertical integration, consider $q$ datasets $\bm{X}^{(k)} \in \mathbb{R}^{n_k \times p_k}$ measuring different variable sets, possibly on non-overlapping subjects. Synthetic data methods aim to learn an estimate $\widehat{\mathcal{P}}$ of the joint distribution over all $p = \sum_{k=1}^q p_k$ variables and generate complete synthetic observations $\bm{\widetilde{X}} \in \mathbb{R}^{n \times p}$. Modern synthetic data methods range from classical techniques based on parametric models to more sophisticated methods leveraging deep generative models such as generative adversarial networks and variational autoencoders. Recent work has also explored the use of synthetic data in semi-supervised learning and prediction-powered inference, where a small labeled dataset $\{(\bm{x}_i, y_i)\}_{i=1}^n$ ($\bm{x}_i \in \mathbb{R}^p$, and $y_i \in \mathbb{R}$) is combined with a larger set of unlabeled data and model predictions $\{\widehat{y}_i\}_{i=1}^N$ to improve estimation and inference \cite{angelopoulos2023prediction,angelopoulos2023ppi++,miao2023assumption,gan2023prediction, miao2024task, zrnic2024active,zrnic2024cross, gu2024local}.

\paragraph{Future directions:}

Current summary statistics are based on strong modeling assumptions, in particular linearity of the effect and no interaction. Open directions include creating novel summary measures or demonstrating how existing summary statistics can be enhanced to model more complex non-linear relationships or interactions between variables. While classical meta-analysis for a single variable of interest provides measures of heterogeneity, measures of between-study variation  become more complex in high-dimensional settings. Nonetheless, when integrating summary-level data across various sources, it is important to account for heterogeneity. Analysis of a high-dimensional set of variables may offer additional opportunities for modeling and learning across covariates. 
% \textcolor{red}{Perhaps this can be linked to Section~\ref{SecHetero}?} 
Systematic missingness or selection bias are common phenomena in health data analytics. When determined by known and measured covariates, these factors can be accounted and adjusted for when individual-level data are available; however, it is impossible to retrospectively adjust summary statistics without access to the initial individual-level data. (See also Section~\ref{SecMissing} on integrating data with missingness.)  
% \textcolor{red}{This might be linked to Section~\ref{SecMissing}.}

Several important questions remain open regarding data integration with synthetic data. First,
characterizing when estimators or test statistics computed on synthetic data $\bm{\widetilde{X}}$
remain consistent or achieve the same rates as those computed on the original data $\bm{X}$ is not
well understood, particularly when the synthetic data distribution $\widehat{\mathcal{P}}$ differs from
the true distribution $\mathcal{P}$. For instance, the use of predicted outcome values for unlabeled data
can result in biased estimation and inference in semi-supervised settings. Some approaches
leverage tools related to classical missing data analysis \cite{robins1995semiparametric} or
prediction-powered inference \cite{angelopoulos2023prediction} to remove the bias while
achieving efficient inference. Other recent proposals have developed semi-parametrically
efficient inference for semi-supervised settings, allowing for a broader range of target
parameters \cite{xu2025unified}. Similar adjustments are needed when the entire dataset used
for estimation and inference, e.g., the set of predictors and outcomes in regression settings,
is synthetically generated. In this fully synthetic setting, naive analysis of synthetic data
is known to yield inconsistent inference with inflated type~1 error, even when point estimates
are unbiased \cite{decruyenaere2023real}. Recent proposals address this by leveraging summary
statistics from the original data or by applying debiasing strategies targeting the specific
analysis of interest to restore standard convergence rates and achieve valid inference
\cite{decruyenaere2024debiasing, keret2025glm}.

A second challenge arises in high-dimensional settings, where determining the sample size
required so that $\widehat{\mathcal{P}}$ preserves sufficient information for downstream tasks such as
estimating a regression parameter $\bm{\beta} \in \mathbb{R}^p$ in a generalized linear model  or conducting a hypothesis test
remains an open problem---especially when integrating across heterogeneous sources with varying
sample sizes $n_j$ and/or noise levels. In fact, prediction-powered inference and inference for
synthetic data have yet to be extended to high-dimensional settings.
Future work could focus on developing formal frameworks for assessing the utility of synthetic
data for specific integration tasks; establishing guarantees for inference procedures based on
integrated synthetic datasets; and designing computationally efficient methods that scale to
high-dimensional multi-source problems.

\subsection{Integrating heterogeneous data}
\label{SecHetero}

The following subsections discuss some common problem settings involving data from heterogeneous distributions, as well as general frameworks for modeling the form of heterogeneity and performing inference on shared parameters. We focus on open directions.

\subsubsection{Transfer learning}
\label{SecTransfer}
\paragraph{Background:}
Transfer learning addresses the problem of leveraging information from source data to improve estimation or prediction on a target domain, particularly when target domain data are limited. Consider $K$ source domains, where the data from the $k$-th domain are given by
\[
\{({\bm x}_i^{(k)}, y_i^{(k)})\}_{i=1}^{n_k} \sim \mathcal{P}_k, \quad k = 1, \ldots, K,
\]
with ${\bm x}_i^{(k)} \in \mathbb{R}^p$ denoting the feature vector and $y_i^{(k)} \in \mathbb{R}$ representing the response. Also, consider a target domain with limited data $\{(\bm{x}_i^{(0)}, y_i^{(0)})\}_{i=1}^{n_0} \sim \mathcal{P}_0$, where $n_0 \ll n_k$. The goal is to estimate a target parameter vector $\bm{\theta}_0$ or prediction function $f_0$ by leveraging information from the source domains. This paradigm has become increasingly relevant in high-dimensional settings where collecting sufficient labeled data for each new task or domain is prohibitively expensive or time-consuming \cite{bastani2021predicting, li2022transfer}. The main challenge in transfer learning is to determine when and how to transfer knowledge between domains, especially when source and target distributions differ. Recent work has focused on characterizing the bias-variance tradeoff inherent in transfer learning, where incorporating information from mismatched source domains can introduce bias, but reduce variance by a factor proportional to the effective sample size $n_0 + \sum_{k=1}^K n_k$ \cite{li2024knowledgetransfer,suder2025bayesian}. Several methods have emerged, including domain adaptation techniques that adjust for distribution shift \cite{reeve2025nonstationary,cai2021transfer}, multi-task learning frameworks that jointly estimate parameters $(\bm{\theta_0}, \bm{\theta_1}, \ldots, \bm{\theta_K)}$ across related tasks \cite{behdin2025multitask,hanneke2022multitask}, and meta-learning techniques that learn transferable representations or algorithms from multiple source tasks \cite{cella2025meta}.

\paragraph{Future directions:}
Recent papers have derived minimax rates for transfer learning in regression, classification, and other supervised learning tasks, characterizing how the amount of source data $\sum_{k=1}^K n_k$, the degree of domain shift---measured in regression settings by the $\ell_q$-norm of the contrast $\max_k \|\bm{\theta_k} - \bm{\theta_0}\|_q$ and in nonparametric classification settings by the strength of the transfer relationship between source and target conditional distributions---and dimensionality $p$ jointly determine the benefit of transfer \cite{li2022transfer, tian2023transfer, cai2021transfer, reeve2021adaptive}. In high-dimensional settings where $p \gg n_0$, transfer learning methods often incorporate sparsity assumptions, where the target parameter $\bm{\theta_0}$ is sparse ($\|\bm{\theta_0}\|_0 \ll p$) and source parameters take the form $\bm{\theta_k} = \bm{\theta_0} + \bm{\delta_k}$, with sparse domain-specific deviations satisfying $\|\bm{\delta_k}\|_q \leq h$ for informative sources \cite{li2022transfer, tian2023transfer}, or low-rank structures where $\bm{\theta_k} \in \mathcal{M}$ for some low-dimensional subspace $\mathcal{M} \subset \mathbb{R}^p$. Computational considerations are also important, as many transfer learning procedures involve solving joint optimization problems of the form $\min_{\bm{\theta_0}, \ldots, \bm{\theta_K}} \sum_{k=0}^K \left(L_k(\bm{\theta_k}) + \lambda R(\bm{\theta_0}, \ldots, \bm{\theta_K})\right)$ over multiple domains or iteratively adapting models from source to target. Open problems include developing robust transfer learning methods that perform well even when some source domains are highly mismatched \cite{li2022transfer}, extending transfer learning theory beyond linear and generalized linear models to modern machine learning architectures, and designing efficient algorithms for large-scale multi-source transfer learning with privacy or communication constraints (cf.\ Section~\ref{SecDistributed}).

\subsubsection{Contaminated data}

\paragraph{Background:}
One setting beyond the usual i.i.d.\ data generating process is Huber's contamination model $(1-\epsilon)\mathcal{P}_{\bm{\theta}}+\epsilon \mathcal{Q}$ \cite{huber1964robust}, where an $\epsilon$ fraction of samples are generated by an arbitrary outlier distribution $\mathcal{Q}$. In recent years, this model has gained much attention in high-dimensional settings where the parameter of interest $\bm{\theta}$ is a location vector, covariance matrix, or regression vector \cite{lai2016agnostic,diakonikolas2019robust,chen2018robust,diakonikolas2023algorithmic}. Both information-theoretic and computational limits have been identified in various settings. From a practical perspective, Huber's model is quite pessimistic in the sense that consistent estimation is impossible when the contamination proportion $\epsilon$ is a constant. In various specific settings, natural conditions on $\mathcal{Q}$ can be formulated, which leads to better estimation and inference guarantees. However, little is known when the problem is high-dimensional.

\paragraph{Future directions:}
Several recent papers consider network analysis with node contamination \cite{acharya2022robust,liu2022minimax}. This leads to a completely new paradigm of robust estimation that is quite different from the sequential setting, where an $\epsilon$ fraction of observations are contaminated, missing, or censored. Network data are naturally organized as an adjacency matrix. If an $\epsilon$ fraction of nodes are contaminated, this means an $\epsilon$ fraction of rows and columns take arbitrary values. In other words, the observed adjacency matrix contains a clean sub-matrix with unknown location. Optimal algorithms for robust recovery in this setting thus need to identify the sub-matrix, leading to completely new challenges. More generally, matrix data with simultaneous row/column contamination can be considered in conjunction with the following problems: ranking from pairwise comparisons \cite{bradley1952rank,gao2023uncertainty}, synchronization \cite{singer2011angular,wang2013exact}, and item response theory \cite{rasch1960probabilistic,nguyen2022spectral}. Another interesting and more general direction is the analysis of tensor data (see Section \ref{ssec:tensor} for more information on tensor data). Analogous to row and column contamination for matrix data, one major challenge is to understand the consequence of mode contamination in tensor data in terms of statistical and computational limits discussed in Section \ref{sec:trade-off}.

%Adaptive inference is a fundamental question in high-dimensional statistics. While adaptive point estimation can usually be solved via standard strategies such as penalized estimation or Lepski's method \cite{lepskii1991problem}, adaptive confidence sets are much more challenging to construct \cite{nickl2013confidence}. In the setting of integrative data analysis, a key question is to quantify the uncertainty of an algorithm with respect to the quality of data, unknown to statisticians. A recent construction of robust confidence intervals for location parameters reveals the necessary cost of adaptation to the unknown contamination level \cite{luo2024adaptive}. Extensions of adaptive confidence intervals to adaptive confidence sets in high dimensions may incur additional computational costs in addition to the cost of adaptation \cite{diakonikolas2023algorithmic}.

\subsubsection{Heteroskedasticity}

\paragraph{Background:}

We consider high-dimensional estimation problems in which samples are drawn from heterogeneous environments but share a common latent structure. A prototypical model is
$
\bm{x}_i \sim \mathcal{N}(\bm{\theta}, \bm{\Sigma}_i), \quad i=1,\ldots,n,
$
where the covariance matrices $\{\bm{\Sigma}_i\}$ are unknown and may vary arbitrarily across samples. The goal is to estimate the common mean $\bm{\theta} \in \mathbb{R}^p$. This model captures a fundamental challenge in modern high-dimensional data analysis: heterogeneity in nuisance parameters across samples. This perspective is related to multiview or multi-environment learning, but differs in that the heterogeneity arises in high-dimensional nuisance parameters rather than multiple structured views. In the one-dimensional setting, this problem reduces to the \emph{entangled Gaussian model} \cite{chierichetti2014learning}, where one observes $x_i \sim \mathcal{N}(\theta, \sigma_i^2)$ with unknown variance profile. A commonly studied special case is the \emph{subset-of-signals} model \cite{liang2020learning}, where an $(1-\epsilon)$ fraction of samples have unit variance and the remaining samples have arbitrary variances. This corresponds to the mixture model $(1-\epsilon)\mathcal{N}(\theta,1) + \epsilon \int \mathcal{N}(\theta,\sigma^2)\, d\Lambda(\sigma)$, which can be viewed as an instance of Huber's contamination model. While sharp statistical guarantees are known in such one-dimensional settings, extensions to high dimensions have only recently begun to be understood \cite{diakonikolas2025entangled}. Beyond Gaussian mean estimation, similar heterogeneity arises in a variety of high-dimensional models.

\begin{enumerate}
\item \textit{Regression with heterogeneous noise \cite{pensia2022estimating}:} $y_i \mid \bm{x}_i \sim \mathcal{N}(\bm{x}_i^T \bm{\beta}, \sigma_i^2)$, where the goal is to estimate $\bm{\beta}$ under unknown variance profile.
\item \textit{Principal component analysis with heterogeneous noise \cite{franks2019shared}:} $\bm{x}_i \sim \mathcal{N}(0, \bm{U}\bm{\Lambda}_i\bm{U}^T + \sigma_i^2 \bm{I}_p)$, where the goal is to estimate the subspace $\bm{U}$.
\item \textit{Models with heterogeneous contamination rates \cite{chaudhuri2025robust}:} $\bm{x}_i \sim (1-\epsilon_i)\mathcal{P}_{\bm{\theta}} + \epsilon_i \mathcal{Q}_i$, extending Huber's model to non-uniform contamination across samples.
\end{enumerate}

\paragraph{Future directions:}
A common theme in the examples above is the presence of sample-wise heterogeneity: observations are generated from different environments but share a common underlying structure. In the Gaussian mean setting, a benchmark estimator with known covariances is the weighted average of the samples. When the covariance profile is unknown, the challenge is to develop procedures that adapt to heterogeneous nuisance parameters while achieving optimal statistical guarantees.

One principled approach to handling such heterogeneity is empirical Bayes (EB). However, the current understanding of EB is largely limited to the one-dimensional Gaussian mean setting \cite{han2026empirical}. Extending EB methodology and theory to high-dimensional regimes and to models beyond mean estimation remains an important open problem. 

On the computational side, it is largely unclear whether polynomial-time algorithms can achieve optimal rates in these heterogeneous high-dimensional settings, particularly when the nuisance parameters are high-dimensional and unstructured. Additional directions include identifying the most informative subset of samples \cite{pensia2022estimating} and developing optimal error criteria tailored to heterogeneous data \cite{compton2024near}.

\subsubsection{Missing data}
\label{SecMissing}

\paragraph{Background:}
A natural connection between missing and contaminated data is that missingness can be regarded as a special instance of contamination. Indeed, if one naively imputes all missing data by some arbitrary value, the missing data can be treated as outliers and algorithms in robust estimation can be applied directly. However, depending on the specific missing mechanism, it is possible that one could do much better for certain estimation tasks \cite{ma2024estimation}.
A related topic is censored data, where one only observes data after performing a truncation mechanism. Some censoring models can be regarded as special cases of missing data models, since if an observation is truncated, one can simply delete it, or equivalently, regard it as missing. However, regarding censored data as missing data can also lead to suboptimal statistical error. It is thus important to understand to what extent better estimation or inference guarantees might be established \cite{bhattacharyya2025learning}.
Important mechanisms of missing/censoring include 1) missing completely at random, 2) missing at random, 3) missing not at random, 4) left and interval censoring, 5) current status data, and 6) self-censoring.

\paragraph{Future directions:}
The above settings can be studied in an analogous way to Huber's model. Consider the setting of data missing not at random. A natural model is that samples are generated from some mixture distribution on $\mathbb{R}^p\cup\{*\}$, such that with $1-\epsilon$ probability, $\bm{x}_i$ generated from $\mathcal{P}_{\bm{\theta}}$ is always revealed; and with $\epsilon$ probability, the sample is missing according to certain probability that may depend on $\bm{x}_i$. One can also define  $\epsilon$ versions of other missing/censoring models. Specifically, problems such as high-dimensional location/scatter estimation, regression and PCA can be considered. In most settings, optimal estimation rates and polynomial-time algorithms are unknown except the recent works of location estimation in a specific setting of missing mechanism \cite{verchand2026high,diakonikolas2026high}. 

%For Ali: add refs/notes about federated learning with distribution shift

%%%%%

\subsection{Distributed learning}
\label{SecDistributed}

\paragraph{Background:}
We now turn to the topic of statistical inference in distributed settings, where data are collected at different locations, with restrictions on the amount of information that may be communicated between locations in a potentially iterative fashion. A primary motivation for this setting is that if a dataset is extremely large, one might not want to store or process the entire dataset at one location (which might correspond to a single computer server). Another benefit of the distributed setting is that locally collected data might contain sensitive information that individuals do not want to share more widely without appropriate postprocessing (in this case, locations might refer to individual hospitals). Finally, a failure of the storage or computational capabilities of a single location would not be nearly as catastrophic as in the case of the usual centralized setting, when one would need to collect additional data and recompute an estimator on the new dataset.

Critical theoretical questions to address include the number of samples required to perform valid statistical inference in a distributed manner (as a function of the number of locations and parameters of the statistical problem) and the amount of communication needed to obtain accurate estimators. In addition, constraints such as privacy or robustness might be placed on information transmitted between local machines; and memory/storage constraints might be placed on local machines or a central server.

\paragraph{Recent work:}
% The survey papers~\cite{gao2022review, chen2014divide, zhou2024distributed} provide an excellent review of existing statistical literature on distributed learning. 
We first provide an overview of several  theoretical questions studied in recent work, with an emphasis on high-dimensional settings; we refer the reader to excellent survey papers~\cite{gao2022review, chen2014divide, zhou2024distributed} for additional references. We also mention the survey paper on federated learning~\cite{li2020federated}, which further provides a computer science perspective on some of these questions. 

\vspace{2mm}
\textit{One-shot averaging}:  Some of the earliest work in the statistics literature concerns ``one-shot averaging," where i.i.d.\ observations are distributed among machines, which compute estimates based on their local observations and then pass these estimates to a central server, which averages the local observations to produce a final estimator. Such an approach was taken in \cite{zhang2013communication, zhang2013information}, with the latter paper imposing an average bit constraint on the amount of communication allowed for the messages passed from local machines to the central server. The paper \cite{rosenblatt2016optimality} further analyzed the tradeoff between the accuracy of such one-shot averages and the complexity of the problem, measured in terms of the number of local machines and the dimensionality of the underlying estimation problem. \cite{minsker2019distributed} studied the benefits of using robust estimators instead of the mean for one-shot aggregation.

\begin{figure}[t]
\begin{center}
\includegraphics[width=7in]{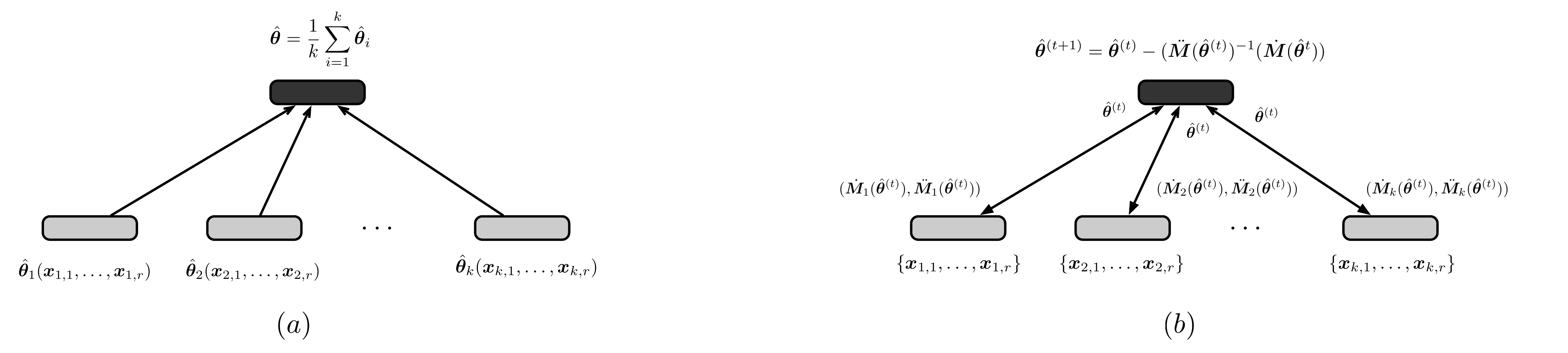}
\caption{A visual comparison of (a) one-shot averaging vs.\ (b) iterative optimization. Machine $j$ observes only the vectors $\{\bm{x}_{j,1}, \dots, \bm{x}_{j,r}\}$. In the case of (b), a new estimate $\bm{\hat{\theta}}^{(t+1)}$ is computed at the central server by averaging gradient/Hessian estimates computed using the local data at each of the machines and then taking a Newton step. Note that the iterates $\{\bm{\hat{\theta}}^{(t)}\}$ must be passed back to all machines at each iteration, in order to compute local gradients and Hessians $(\bm{\dot{M}}_i(\bm{\hat{\theta}}^{(t)}), \bm{\ddot{M}}_i(\bm{\hat{\theta}}^{(t)}))$, which are then aggregated, e.g., using simple averaging, to obtain $(\bm{\dot{M}}(\bm{\hat{\theta}}^{(t)}), \bm{\ddot{M}}(\bm{\hat{\theta}}^{(t)}))$.}
%\textcolor{red}{If possible, can we make $\theta$s in the figure boldfaced? Also it would be great if we specify in the caption how we would like to see $X_{ij}$. It seems to be me that these are modeled as vectors? If so, can we make them lowercase-boldfaced, e.g. $\bm{x}_{1,1}$? If you would like to see them as  matrices, then the notation is fine. }}
\end{center}
\end{figure}

\vspace{2mm}
\textit{Iterative optimization}: When information may be exchanged multiple times between local machines and the central server, \cite{huang2019distributed} proposed and studied a one-step estimator, where an initial estimator $\bm{\widebar{\theta}}$ of the parameters of interest is first computed based on one-shot averaging, then passed to local machines, which compute estimates of the Hessian and gradient of the loss function at $\bm{\widebar{\theta}}$ based on their local samples. These estimates are then passed back to the central server, which computes an average Hessian and gradient and takes one Newton step.

In order to reduce computation and communication costs, \cite{jordan2019communication} proposed a method which only involves computing a Hessian matrix at one local machine based on the initial estimator $\bm{\widebar{\theta}}$, which is then passed to the central server. They analyzed the convergence rate of the estimator as a function of the number of approximate Newton steps taken. \cite{fan2023communication} incorporated regularization into the update steps in order to relax the conditions required for convergence, and extend the applicability of the iterative optimization method to higher dimensions. 

\vspace{2mm}
\textit{High-dimensional regression}: The papers \cite{lee2017communication, battey2018distributed} studied the problem of distributed high-dimensional linear regression, where $n$ i.i.d.\ observations from a sparse linear model are divided among $k$ machines. For testing hypotheses concerning the values of certain coordinates, they proposed a test based on averaging standardized debiased Lasso statistics computed on local machines. They also studied the statistical error of the regression estimator taken as an average of Lasso estimators computed on local machines. \cite{lee2017communication} further noted that a debiased Lasso estimator may be computed using only local debiasing steps that adjust a subset of the coordinates of the central estimator based on observations on each local machine, thus further reducing the computational and communication complexity of the overall algorithm. They presented a communication-efficient procedure when two rounds of communication are permitted between local machines and the central server. The paper~\cite{song2015split} proposed a Bayesian method for high-dimensional linear regression under an alternative scenario where data are split between machines by partitioning the coordinates rather than the samples.

\vspace{2mm}
\textit{PCA}: The paper~\cite{fan2019distributed} studied the problem of extracting top principal eigenvectors of a dataset, where the i.i.d.\ observations are divided among local machines. Their proposed estimator consists of averaging a low-rank approximation of the sample covariance matrices computed on local machines, and then extracting the top eigenvectors of the aggregated covariance matrix. The paper discussed a few other approaches to distributed PCA which have appeared in recent literature, including methods that analyze data which are partitioned by features rather than observations (vertical vs.\ horizontal partitioning), but pointed out that few of these papers study the statistical error of the proposed methods. 

\vspace{2mm}
\textit{Quantile regression}:  Distributed quantile regression has received special attention, due to the mildly non-smooth nature of the loss function that causes complications in optimization-based approaches. \cite{tan2022communication} studied a method based on smoothing the objective function. \cite{yang2025quantile} studied high-dimensional quantile regression using a smoothed quantile loss and a concave penalty. \cite{chen2019quantile} studied an interesting version of quantile regression which involves one central machine with a \emph{memory} constraint, meaning it can only store a limited batch of data at any time. They proposed a method based on distributing the data into smaller chunks, then iteratively updating the estimator as successive chunks of data are processed. 

\vspace{2mm}
\textit{Sign selection}: The paper~\cite{liu2024majority} studied the problem of estimating the sign pattern of nonzero components of a high-dimensional parameter vector when i.i.d.\ data are distributed among local machines. In addition, they leveraged the distributed nature of the problem to develop a version of a distributed sign selection algorithm that also satisfies a notion of group privacy. The algorithm is based on a majority vote of local sign estimates obtained via regularization, combined with a private peeling algorithm which adds Laplace noise to vote counts in order to ensure privacy.
The same technique of majority voting appeared in the earlier paper~\cite{chen2014split}, which studied distributed methods for sign selection in high-dimensional linear models, without a privacy constraint.

\vspace{2mm}
\textit{Nonparametric estimation}: 
\cite{zhu2018distributed, szabo2020adaptive, cai2022distributed} studied the problem of nonparametric regression of a smooth function, when i.i.d.\ data are distributed among local machines. The number of bits communicated from each local machine to the central server in a one-shot manner must satisfy a certain bound (required to hold in expectation in the papers \cite{szabo2020adaptive, cai2022distributed}). The papers derived minimax rates of estimation in terms of the number of samples, number of machines, communication constraints, and smoothness of the function class. \cite{szabo2020adaptive, cai2022distributed} also characterized the cost of adaptation, when the level of smoothness of the true regression function is unknown. 

\vspace{2mm}
\textit{Hypothesis testing}:  Distributed hypothesis testing is studied in the papers~\cite{szabo2022optimal, szabo2023optimal}, which focus on signal detection in high-dimensional and nonparametric settings. They studied the minimax testing error under a bit constraint on the messages communicated to the central server in a single iteration. The paper~\cite{vuursteen2024optimal} studied communication-constrained goodness-of-fit testing for discrete distributions, while the paper~\cite{cai2024federated} considered privacy-constrained testing for nonparametric distributions. In a separate line of work, the paper~\cite{vuursteen2023optimal} studied optimal methods for aggregating test statistics from locally computed hypothesis tests based on disparate batches of data, where the communication constraint is imposed in the form of the expected absolute value of individual test statistics satisfying a prescribed bound.

\cite{pensia2023communication} studied the question of hypothesis testing between discrete distributions, where i.i.d.\ observations are observed after being passed through communication-constrained channels. This framework, which is a special case of distributed inference where each local machine only observes a single sample, has also received a fair amount of attention in the information theory literature~\cite{acharya2020inference, han2018geometric}.

\paragraph{Future directions:}
We now outline some future research directions below. 

\vspace{2mm}
\textit{More refined tradeoffs}: While a few of the aforementioned papers derive both upper and lower bounds on the performance of inference procedures corresponding to the problem settings, a fair amount of existing work focuses on presenting algorithms for distributed inference and deriving accuracy bounds, leaving open the question of whether other algorithms exist with better performance. Thus, an interesting future direction that would likely lead to the development of new lower bound techniques and new statistical methodology would be to revisit the papers above and study the question of optimality. This would, in turn, provide a more precise characterization of the tradeoffs in statistical performance caused by restricting attention to decentralized procedures. In particular, it would be important to have a more rigorous understanding of the amount of communication used by algorithms (e.g., in terms of bits), rather than simply counting iterations or dimensionality.

\vspace{2mm}
\textit{More complex communication protocols}: Taking a cue from the information theory literature, one might try to determine optimal procedures according to the number of bits required in communication, or imposing other constraints such as sparsity of message vectors~\cite{wang2018atomo}. Pushing this idea even further, one could require information to be transmitted in a fault-tolerant manner, which is an active area of research in federated learning~\cite{li2018review}, where it is important for the system to perform well even when a number of local machines fail to report their estimates at a given iteration. From a more traditional statistical perspective, one might treat these failed estimates as missing data, or in the sense of statistical robustness~\cite{minsker2019distributed, van2025robust}, but a coding-theoretic perspective for communication-constrained inference might be interesting to explore~\cite{zhang2024coded}. 

\vspace{2mm}
\textit{Sequential problems}: In many contemporary applications, data are acquired in a temporal manner. This opens the door to more sophisticated estimation procedures than have been studied in most previous literature, which assumes that i.i.d.\ draws are collected and distributed to local machines all at once. In particular, if data were collected sequentially, local nodes could make decisions about whether to collect additional data points or change the type of information transmitted to the central server, based on information they have received from other machines and/or the central server.
Such questions are related to online learning~\cite{shamir2014fundamental}, the role of interactivity in communication-constrained inference~\cite{acharya2022role}, and memory-constrained inference~\cite{steinhardt2015minimax, dudeja2024statistical}, but existing work only scratches at the surface of what sorts of  settings might be studied.

In a more classical statistical framework, a natural direction would be to develop methods for detecting changes in the underlying data generating process when i.i.d.\ data are distributed among local machines. This area seems to have received very little attention~\cite{yang2024communication}. Other questions such as time-series segmentation, forecasting, or classification could also be posed in a distributed setting. 

\begin{figure}[t]
\begin{center}
\includegraphics[width=6in]{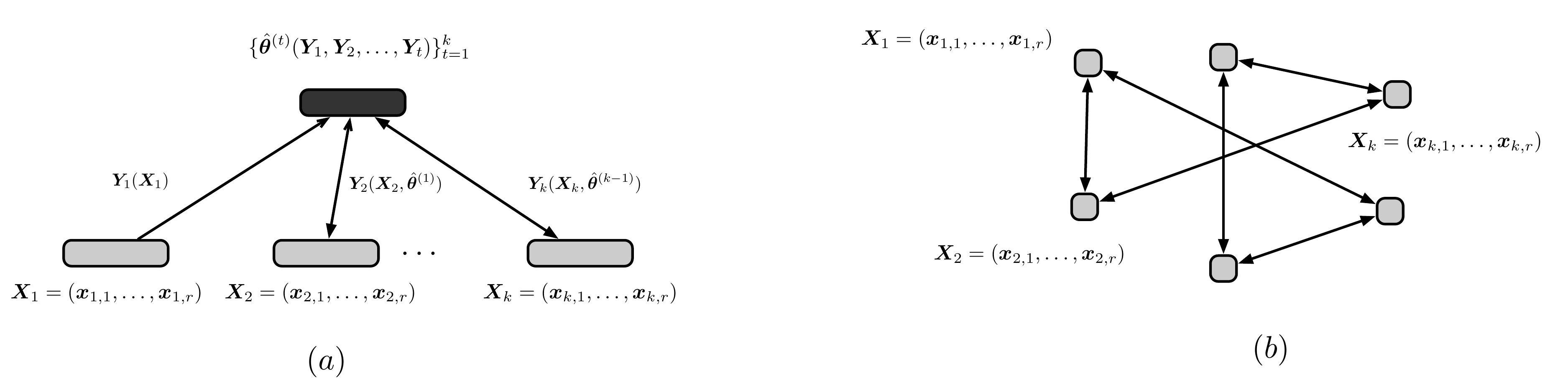}
\caption{Some interesting open directions in distributed learning. (a) An illustration of a sequential setting, where the messages $\{\bm{Y}_t\}$ (which may contain information such as gradients, Hessians, etc.) sent from local machines to the central server can depend on their own observations as well as previously released estimates $\{\bm{\theta}^{(t-1)}\}$ transmitted from the central server, calculated using previously received messages from local machines. Local machines could also choose to adapt their data collection strategy based on past messages. (b) An illustration of a fully decentralized setting, which does not involve communication to/from a central server.}
\end{center}
\end{figure}

\vspace{2mm}
\textit{Bayesian distributed learning}: Relatively little work has been done in a Bayesian framework for distributed inference, and distributed learning is recognized as an important topic in Bayesian statistics with many unanswered questions~\cite{winter2024emerging}. A natural question is how to compute an approximate posterior from i.i.d.\ data when the data are stored on local machines and only local computations (e.g., calculating or sampling from a posterior based on local data) can be performed.

One line of work~\cite{minsker2017robust} proposes and studies aggregation procedures for ``subset posteriors" computed using data on local machines. In one approach, the same prior is used for all local posterior calculations, while in another approach, the prior $\pi(\bm{\theta})$ is incorporated as a ``fractional prior" of the form $\pi(\bm{\theta})^{1/k}$ among the $k$ machines. Aggregation procedures based on the mean and median (in a suitable Wasserstein metric) have been studied. However, theoretical results for the convergence of the resulting posterior are quite limited, as they depend on normality of the subset posteriors, and the methods are not readily scalable to high dimensions. Another approach, based on stochastic gradient Markov Chain Monte Carlo, proposes to sample data from the full posterior on all samples by computing and communicating gradients on separate batches of data~\cite{ahn2014distributed}. However, as such a procedure would involve multiple rounds of communication between local machines and the central server, its practicality remains fairly limited.

Evidently, many opportunities exist for further theoretical development in the field of Bayesian distributed computation. Extending the scope of existing theory for subset posterior averaging or rigorously characterizing the amount of communication needed for stochastic gradient methods in low dimensions would already constitute an important advancement, whereas devising methods that are practical for high-dimensional problems could prove significantly more complex and challenging. 

\vspace{2mm}
\textit{Different network topologies}:  
All of the work mentioned above assumes the existence of a central server, which interacts with local machines (possibly in an iterated manner). \emph{Decentralized learning} has been studied elsewhere in engineering and computer science, and concerns a different paradigm where local machines communicate over the edges of a predetermined network, and maintain local estimates at each step which are functions of their previous estimates and those of their neighbors. The most commonly studied problem, popularized by work in federated learning, is decentralized optimization~\cite{cevher2014convex}, where the goal is to minimize the sum $L(\bm{\theta}) = \sum_{k=1}^K L_k(\bm{\theta})$ of local loss functions, where each of the $K$ nodes in the network sees a separate batch of i.i.d.\ data from an unknown distribution $\mathcal{P}_k$. Such a paradigm of decentralized learning has seen very little attention in statistics, save the paper~\cite{gu2024statistical}, which studies statistical properties of the iterates of a decentralized version of stochastic gradient descent, with local estimates obtained by averaging gradient directions among neighboring nodes.

Moving beyond this setting, many open problems remain concerning how to perform estimation for different statistical problems in a decentralized manner (not necessarily via minimizing a sum of local loss functions), and what other sorts of information might be communicated between nodes to efficiently solve a joint inference problem. A different branch of research in the engineering community is the study of \emph{gossip algorithms}, where some theoretical results have been derived concerning the estimation of certain statistical problems relevant to sensor network analysis based on local communication~\cite{dimakis2010gossip}. This could be a fruitful area from which one could extract additional paradigms of communication and distributed information that would be interesting to study from a modern statistical learning perspective.
%\documentclass{article}

% Language setting
% Replace `english' with e.g. `spanish' to change the document language
%\usepackage[english]{babel}

% Set page size and margins
% Replace `letterpaper' with `a4paper' for UK/EU standard size
%\usepackage[letterpaper,top=2cm,bottom=2cm,left=3cm,right=3cm,marginparwidth=1.75cm]{geometry}

% Useful packages
%\usepackage{amsmath, amsfonts}
%\usepackage{bm}
%\usepackage{graphicx}
%\usepackage[colorlinks=true, allcolors=blue]{hyperref}

%\newtheorem{example}{Example}

%Notations
%\newcommand{\bx}{\bm{x}}
%\newcommand{\bbeta}{\bm{\beta}}
%\newcommand{\mX}{X}
%\newcommand{\bv}{\bm{v}}
%\newcommand{\bu}{\bm{u}}
%\newcommand{\mZ}{Z}
%\newcommand{\R}{\mathbb{R}}

\section{Challenges in high-dimensional asymptotic analysis}
\label{sec:asymptotics}
%\author{Sivaraman Balakrishnan\thanks{Department of Statistics \& Data Science and Department of Machine Learning, Carnegie Mellon University},~ Rishabh Dudeja\thanks{Department of Statistics, University of Wisconsin-Madison},~ Chao Gao\thanks{Department of Statistics, University of Chicago},~ Jason M. Klusowski\thanks{Department of Operations Research and Financial Engineering, Princeton University},~ \\ Po-Ling Loh\thanks{Statistical Laboratory, Department of Pure Mathematics and Mathematical Statistics, University of Cambridge}, ~ Arian Maleki\thanks{Department of Statistics, Columbia University}, ~ Subhabrata Sen\thanks{Department of Statistics, Harvard University},~ Ali Shoajaie\thanks{Department of Biostatistics and Department of Statistics, University of Washington, Seattle}, ~ \\
%Weijie Su\thanks{Department of Statistics \& Data Science, Wharton School, University of Pennsylvania}, ~ Anru Zhang\thanks{Department of Biostatistics \& Bioinformatics and Department of Computer Science, Duke University},~ Verena Zuber\thanks{Department of Epidemiology and Biostatistics, School of Public Health, Imperial College London}}

\subsection{Introduction}

 Asymptotic analysis has long been a cornerstone of classical statistical theory. Its strength lies in providing general and tractable approximations for the behavior of estimators, such as their risk or distribution, and for the error rates of tests, particularly when exact finite-sample results are unavailable or intractable. This allows us to compare competing estimators/tests, and determine the optimal procedure in a given class. 
 A prominent example is maximum likelihood estimation (MLE), one of the most widely used methods in classical statistics, whose theoretical justification rests largely on asymptotic results: under suitable regularity conditions, the MLE is consistent, asymptotically normal, and asymptotically efficient.

The success of asymptotic analysis in classical statistics has motivated researchers to investigate extensions to high-dimensional problems. However, high-dimensional settings demand a more delicate and nuanced asymptotic treatment of estimators and tests. Key challenges include the growth of the number of parameters relative to the sample size, the intrinsic complexity of the parameter space, and the need for regularization. To illustrate some of these complications and set the stage for the challenges that arise in this context, we review two canonical examples that have been extensively studied in the literature: (ii) additive spiked models and (ii) sparse linear regression.

 \subsubsection{Principal component analysis}
 \label{sec:pcaexample}
 One of the simplest models for studying principal component analysis is the \emph{additive spiked model}, in which the goal is to estimate the unit-length $p$-dimensional signal direction $\bv_*$ from an observation of the form
\[
\mX = \lambda \bu \bv_*^{\top} + \mZ,
\]
where $\mX \in \R^{n \times p}$ is the data matrix, 
$\mZ_{ij} \overset{\mathrm{i.i.d.}}{\sim} \mathcal{N}(0,\frac{1}{n})$ denotes additive Gaussian noise, 
and $\bu \in \R^{n}$ has entries $\bu_i\overset{iid}{\sim} \mathcal{N}\left(0,\,\frac{1}{n}\right)$. The parameter $\lambda>0$ denotes the SNR in the problem; the larger $\lambda$, the easier it is to recover the latent signal $\bv_*$. The vector $\bu$ is a latent nuisance parameter. High-dimensional asymptotic approximations are particularly relevant when the number of samples $n$ and the dimension of the parameter space $p$ are large.

At this point, one can ask two different questions: (i) What is the performance of an estimator of interest? For example, 
one can estimate $\bv_*$ using the principal component of $\mX^{\top}\mX$, denoted by $\hat{\bv}$, and then ask how accurate this estimate is. (ii) What is the optimal estimator of $\bv_*$? On a related note, one can also investigate optimal estimators in a sub-class, e.g., the class of low-degree polynomials (see \cite{montanari2025equivalence} and the references therein). 

To study these questions under an asymptotic framework, one needs an appropriate notion of loss. Prominent examples include:
\begin{itemize}
\item \emph{Weak recovery}: Is there an estimator $\hat{\bv}$ that achieves higher correlation than a purely random guess?
\item \emph{Near exact recovery}: Is there an estimator $\hat{\bv}$ such that $\mathcal{L}(\bv_*, \hat{\bv}) \to 0$, where $\mathcal{L}$ represents a generic loss?
\item \emph{Exact recovery}: Is there an estimator $\hat{\bv}$ such that $\mathbb{P}(\hat{\bv} = \bv_*) \to 1$ as $n,p \to \infty$?
\end{itemize}

The questions phrased above can be tackled using appropriate asymptotic frameworks. 
%Different asymptotic frameworks can be employed to analyze the accuracy of this estimator. 
To specify an asymptotic regime, one needs to determine how the dimension $p_n$ and the signal strength $\lambda_n$ scale as $n \to \infty$.
As an example, we collect a concrete statement on  the performance of the principal component of $\mX^{\top}\mX$ under one such asymptotic regime \cite{10.1214/009117905000000233,benaych2011eigenvalues}:

\begin{example}
\label{ExHDA}
Suppose $\frac{p}{n} \to \delta \geq 0$ as $n \to \infty$, while $\lambda_n = \lambda$ remains constant. If $\lambda > {\delta}^{\frac{1}{4}}$, we have
\[
|\langle \hat{\bv}, \bv_* \rangle |
\xrightarrow{p} \sqrt{
\frac{\,1 - \frac{\delta}{\lambda^{^4}}\,}{\,1 + \frac{\delta}{\lambda^{2}}\,}}.
\]
In contrast, if $\lambda < \delta^{\frac{1}{4}}$, we have
\[
\langle \hat{\bv}, \bv_* \rangle \xrightarrow{p} 0 .
\]
\end{example}

\subsubsection{Sparse linear regression}\label{ssec:linear} 
In Example~\ref{ExHDA}, we did not impose any additional structure on $\bv_*$, such as sparsity. 
As a result, the range of possible asymptotic regimes was rather limited. 
We now consider another canonical example in which the asymptotic analysis is richer.

Consider a dataset $\{(y_{1}, \bx_{1}), (y_{2}, \bx_{2}), \ldots, (y_{n}, \bx_{n})\}$, 
where $y_{i} \in \R$ is the response variable and 
$\bx_{i} \in \R^{p}$ is the feature vector of the $i^{\text{th}}$ observation. Suppose the data points are i.i.d.\ and follow the linear model
\begin{align}
y_{i} = \bx_{i}^{\top} \bbeta_* + z_{i}, \label{eq:regression}
\end{align}
where $z_{i} \overset{\text{i.i.d.}}{\sim} \mathcal{N}(0, \sigma^{2})$ is Gaussian noise 
and $\bbeta_* \in \R^{p}$ denotes the unknown parameters to be estimated.
 In \emph{sparse linear regression}, it is assumed that most of the components of $\bbeta_*$ are zero, with only $k$ nonzero coefficients. 
Classical asymptotic theory typically treats $p, k,$ and $\sigma$ as fixed, while $n \to \infty$. However, this regime is not suitable for high-dimensional settings; instead, it is natural to let $p, k,$ and possibly $\sigma$ also grow with $n$. A theoretical statistician faces a range of choices:
\begin{enumerate}
\item How do the quantities $p_{n}$, $k_{n}$, and $\sigma_{n}$ scale with $n$?

\item Given that $p_{n}$ grows with $n$, what assumptions should be made about the feature vectors $\bx_{i}$? 
For instance, one may specify how quickly $\|\bx_{i}\|_{2}^{2}$ grows with $n$, 
or impose additional structural or distributional assumptions on $\bx_{i}$, such as sub-Gaussianity.
\end{enumerate}

We now discuss two of the simplest asymptotic results in sparse linear regression:

\begin{example}
[Consistency regime] 
\label{slr:consistency:1}
Consider the case where $k_{n}, p_{n} \to \infty$ as $n \to \infty$, while $\sigma_{n} = \sigma$ remains fixed. 
Further assume that $\bx_{i} \sim N\!\left(0, \frac{1}{n} \bm{I}_{p}\right)$, and $k/n \to 0$ and $k \log p \,/\, n \to 0$. 
Then
\[
\frac{\inf_{\hat{\bbeta}} \; \sup_{\bbeta} \; \mathbb{E}\bigl[\|\hat{\bbeta} - \bbeta\|_{2}^{2}\bigr]}
     {2 k \log \frac{p}{k}}
\;\longrightarrow\; 1 .
\]
In this consistency regime, the minimax error of $2 k\log\frac{p}{k}$ can be achieved asymptotically via the SLOPE method \cite{bogdan2015slope}, which finds the estimate by solving
\[
\arg\min_{\bbeta} \left\{\frac{1}{2}\sum_{i=1}^n (y_i - \bx_i^T \bbeta)^2 + \sum_{j=1}^p \lambda_j |\bbeta|_{(j)}\right\},
\]
where $\lambda_1 \ge \lambda_2 \cdots \ge \lambda_p \ge 0$ is a non-increasing regularization sequence and $|\bbeta|_{(j)}$'s are the order statistics of $\bbeta$ in absolute value. The optimal risk can be achieved by taking a certain regularization sequence that is inspired by the Benjamini--Hochberg procedure \cite{su2016slope}.
\end{example}
A proof of the above result, along with a review of the existing literature on the topic, can be found in \cite{guo2024note}. Our next result focuses on the LASSO algorithm for estimating ${\bbeta}$:
\[
\hat{\bbeta}^L\in \arg\min_{\bbeta} \left\{\frac{1}{2}\sum_{i=1}^n (y_i - \bx_i^T \bbeta)^2 + \lambda \|\bbeta\|_1\right\}.
\]

\begin{example}
[Inconsistency regime for LASSO \cite{bayati2011dynamics,donoho2011noise}]
\label{example:LASSO}
Consider the case where $k_n, p_n \rightarrow \infty$ as $n \rightarrow \infty$, while $\sigma_n = \sigma$ remains fixed. Further assume that $\bx_i \sim \mathcal{N}(0, \frac{1}{n}\bm{I}_p)$, and $k/p \rightarrow \epsilon$, while $n/p \rightarrow \delta$. For each dimension $p_n$ let the empirical cdf of $\bbeta \in \R^{p_n}$ converge to the probability measure $p_{B}(b)$ weakly. Further suppose that
\[
\frac{1}{p} \|\bbeta\|_2^2 \rightarrow \mathbb{E}(B^2),
\]
where $B\sim p_{B}(b)$. Then
\[
\frac{1}{p}\|\hat{\bbeta} -\bbeta\|_2^2 \overset{p}{\rightarrow} \mathbb{E} (\eta (B+ m Z ; \tau m) - B)^2,
\]
where $\eta (u; \tau m) = (|u|- \tau m)_+ {\rm sign} (u)$, and $\tau$ and $m$ satisfy the following fixed-point equations:
 \begin{eqnarray}
 m = \frac{1}{\delta} \mathbb{E} (\eta(B+m Z; \tau m)- B)^2, \nonumber \\
 \lambda = \tau m (1 - \mathbb{E} (\eta'(B+mZ; \tau m))). 
 \end{eqnarray}
\end{example}

To characterize the performance of statistical procedures in high-dimensions, researchers employ diverse tools arising in high-dimensional probability, information theory, and related disciplines. In particular, recent advances in high-dimensional asymptotics rely heavily on tools such as concentration inequalities, random matrix theory, Gaussian process theory, and empirical process theory. At the same time, questions arising from high-dimensional learning often directly motivate the development of new tools in these areas. 

%In Section~\ref{ssec:hd:asym:tool}, we review some existing tools for high-dimensional analysis, and discuss possible extensions to broaden the applicability of these tools. Subsequently, we highlight 

In this section, we focus on the following open directions: In Section \ref{ssec:hd:asym:tool}, we review some existing tools for high-dimensional asymptotic analysis and discuss challenges associated with using them. In Section \ref{ssec:bayesian}, we discuss questions in high-dimensional Bayesian asymptotics, which are critical in developing a notion of statistical optimality in high-dimensional problems. In Section \ref{ssec:Bayesian:functional}, we discuss the problem of estimating low-dimensional functionals of the regression coefficients under proportional high-dimensional asymptotics.
In Section \ref{sec:designs_questions}, we highlight some research questions seeking to extend the scope of high-dimensional asymptotic guarantees, so that they can be applied out-of-the-box. Finally, in Section~\ref{ssec:refinedapprox}, we discuss opportunities to refine the current high-dimensional asymptotic framework, allowing one to compare across a wider class of competing methods and derive a finer understanding of their statistical performance. 

\subsection{Tools for High-Dimensional Asymptotics} \label{ssec:hd:asym:tool}

\paragraph{Background:}
The markedly different behavior of high-dimensional asymptotics, compared with classical low-dimensional regimes, has motivated the development of new analytical tools. Many of these tools draw inspiration from ideas in statistical physics. In particular, the replica method provides a powerful, though non-rigorous, approach for deriving sharp asymptotic risk characterizations for M-estimators, as well as Bayesian estimators such as the posterior mean. This technique is highly general and can yield predictions even when the M-estimator is defined through a nonconvex loss, or when the prior or likelihood is misspecified in the Bayesian setting. We refer the reader to \cite{nishimori2001statistical,bandeira2018notes,montanari2024friendly} for accessible tutorial introductions to the replica method.

In recent years, several techniques have been developed to rigorously justify predictions inspired by the replica method in certain settings. These include (i) approximate message passing, (ii) Gaussian comparison inequalities, (iii) leave-one-out (or cavity) methods, and (iv) interpolation techniques.
The goal of this section is to review these tools and discuss their  limitations.

\paragraph{Recent work and future directions:}
Below, we review some of the existing tools for high-dimensional asymptotic analysis and discuss a few open challenges in each area.

\begin{itemize}
    \item \emph{Approximate Message Passing (AMP):} AMP algorithms \cite{donoho2009message} are iterative algorithms that are closely related to gradient descent and proximal gradient descent, which are commonly used to compute M-estimators. AMP updates differ from gradient descent updates by a carefully designed correction term called the \emph{Onsager correction}, which ensures that the AMP iterates can be approximately decomposed into a signal term in the direction of the unknown parameter and a Gaussian noise term with i.i.d. coordinates \cite{donoho2009message,bolthausen2009high,bayati2011dynamics}. This \emph{state evolution} decomposition enables an exact characterization of the risk of AMP iterates in the high-dimensional limit. AMP algorithms are often used as a proof technique to derive risk characterizations for M-estimators that minimize a convex loss. The basic proof idea, proposed in \cite{bayati2011dynamics}, is to write down an AMP algorithm that converges to the global minimizer of the loss function. This step often requires the loss to satisfy strong convexity or restricted strong convexity. The risk of the M-estimator can then be derived from the state evolution of the AMP algorithm. We refer the reader to \cite{feng2022unifying} for a tutorial introduction to AMP  and its analysis. 
    The behavior of AMP algorithms for non-convex problems is not well-understood, and investigating this is an important direction for future work.  In some non-convex problems, the iterates of a suitably designed AMP algorithm converge. These problems are called \emph{replica symmetric}, and in these cases, the AMP state evolution gives the correct characterization of the risk of the M-estimator \cite{vilucchio2025asymptotics}. In other non-convex situations, AMP iterates may fail to converge. Is it possible to design an alternative algorithm that finds the global minimizer of a non-convex loss in these situations and use it to analyze the behavior of the corresponding M-estimator? 
    \item \emph{Gaussian Comparison Inequalities:} A different approach to derive risk characterizations for convex M-estimators is called the Convex Gaussian Min-Max (CGMT) framework and is based on Gordon’s comparison inequality \cite{gordon1988milman}. Building on earlier  analyses of phase 
transitions in $\ell_1$ minimization and random convex programs 
\cite{stojnic2009various,oymak2010new,chandrasekaran2012convex,amelunxen2014living}, this framework was developed over a sequence of works \cite{stojnic2013framework,stojnic2013regularly,stojnic2013upper,oymak2013squared,thrampoulidis2015regularized,thrampoulidis2018precise,thrampoulidis2016recovering}. Considering the linear model \eqref{eq:regression} with $n \times p$ i.i.d.\ Gaussian design matrix $\mX$ as an example, the framework shows that the asymptotic properties (such as the risk) of the M-estimator obtained by minimizing a \emph{convex} loss such as the regularized least squares loss and a \emph{convex} regularizer $r: \R^p \rightarrow \R$:
\begin{align} \label{eq:rls}
\hat{\bbeta} & =  \arg\min_{\bbeta \in \R^p} \left\{\frac{1}{2} \| \mathbf{y} - \mX \bbeta\|_2^2 + \lambda \,  r(\bbeta) \right\}
%\sum_{i=1}^p r(\beta_i)
\end{align}
are close to the asymptotic properties of the minimizer of a much simpler optimization problem, where the $n \times p$ design matrix $\mX$ has been replaced by two independent Gaussian random vectors $\bm{g}  \sim \mathcal{N}(\bm{0},\bm{I}_n)$ and $\bm{h} \sim \mathcal{N}(\bm{0}, \bm{I}_m)$. The latter can be reduced  to a 
low-dimensional optimization problem over only a few scalar variables, which can often be solved 
explicitly or, more generally, characterized implicitly via first-order stationarity conditions, ultimately yielding 
exact risk and distributional characterizations of the M-estimator.  The framework crucially relies on the convexity of the objective~\eqref{eq:rls}. For non-convex problems, Gordon's inequality provides a lower bound on the minimum value of the loss function. For some non-convex problems, this lower bound is in fact tight, and one can combine Gordon's inequality with the AMP framework to obtain a matching upper bound and characterize the risk of the non-convex M-estimator \cite{vilucchio2025asymptotics}. For more complicated non-convex problems, the lower bound given by Gordon's inequality is not expected to be tight \cite{zadik2019improved}. 
%
%For non-convex problems, Gordon's inequality provides a lower bound on the minimum value of the loss function. For some non-convex problems this lower bound is tight, and one can combine Gordon's inequality with the AMP framework to obtain a matching upper bound and thereby characterize the risk of the non-convex M-estimator \cite{vilucchio2025asymptotics}. For more complicated non-convex problems, however, the lower bound from Gordon's inequality is not expected to be tight.
In a series of recent works, Stojnic \cite{stojnic2016fully,stojnic2023fully,stojnic2023fully2} developed a clever lifting technique that applies a Laplace-transform-like trick to Gordon's original inequality, yielding tighter lower bounds on several non-convex instances. An important open question is whether this or a similar approach can be made as transparent and widely-accessible as the CGMT, lowering the barrier to adoption by the broader community.
It is worth noting that even within the convex setting, a major challenge remains: extending the CGMT to the multivariate, matrix-valued analogue of \eqref{eq:rls}, where the regression coefficient becomes a $p \times k$ matrix $\bm{B} = [\bbeta_1, \ldots, \bbeta_k]$ and the response becomes an $n \times k$ matrix $\bm{Y} = [\mathbf{y}_1, \ldots, \mathbf{y}_k]$. This matrix generalization, under general loss functions, subsumes multivariate regression, multiclass classification, and the graphical LASSO as special cases. Only partial results are currently available \cite{thrampoulidis2020theoretical}. Can a general, versatile matrix analog of the CGMT be developed to cover these cases?
% Can we develop a similarly general and versatile framework, such as the CGMT framework, to handle general non-convex problems? 

\item \emph{Leave-One Out or Cavity Method:} This technique, which originated in statistical physics and was introduced in statistics by \cite{el2013robust,karoui2013asymptotic,el2018impact,sur2019modern,chen2021spectral, rad2020scalable}, not only provides asymptotic risk characterizations for convex M-estimators, but also characterizes the distribution of individual coordinates, which can be used for inference. Using the regularized least squares (RLS) estimator \eqref{eq:rls} as an example, the leave-one-out method approximates the distribution of the $i^{\text{th}}$ coordinate $\hat{\bbeta}_i$ of the RLS estimator by comparing the full estimator to the leave-one-feature-out estimator $\hat{\bbeta}_{-i}$, obtained by removing feature $i$. A Taylor expansion of the loss around $\hat{\bbeta}_{-i}$ gives an explicit approximation for $\hat{\bbeta}_i$ in terms of the $i^{\text{th}}$ feature vector $\bm{X}_i$ and $\hat\bbeta_{-i}$. Because $\bm{X}_i$ and $\hat\bbeta_{-i}$ are independent, this approximation implies that the distribution of $\hat{\bbeta}_i$ is determined by the empirical distribution of residuals from the leave-one-feature-out problem. Assuming the optimization problem is stable, these residuals are close to those of the original regression problem, which can be analyzed via a leave-one-observation-out argument (dropping observation $j\in [n]$). Together, these steps yield a closed system of fixed-point equations that determine the empirical distributions of the RLS estimator and its residuals, along with their coordinate-wise distributions. Convexity is essential to ensure that the optimization problem is stable under feature/observation removal and the fixed-point equations are closed and uniquely solvable. Versions of the leave-one-out method have also been developed to analyze posterior distributions (with possible misspecification), provided the posterior is log-concave \cite{barbier2025performance,saenz2025characterizing}. Extending the leave-one-out method beyond convex losses and log-concave posteriors is an important direction for future work. 
\item \emph{Interpolation Method:} The interpolation method, developed over a series of works (see e.g.\cite{barbier2019adaptive,barbier2019optimal}), is a general approach to characterize the risk of the posterior mean for high-dimensional Bayesian problems. The method has roots in the study of mean-field spin glass models. Formally, one compares the log-normalization constant of the posterior distribution to that of a simpler distribution (often a product distribution). If the difference between the two normalization constants is small, one can often read off the estimation risk in the posterior of interest from that of the simpler model. The current version of this technique works when both the prior and  likelihood are correctly specified. We refer the interested reader to \cite{montanari2024friendly}  for illustrations of this method. %However, the replica method provides predictions for the risk of the Bayes estimator even when the prior or the likelihood is misspecified. Proving these predictions is an important open problem.
Extending the interpolation method to misspecified problems is an important direction for future inquiry. See \cite{camilli2022inference,guionnet2025estimating} for some recent results in this direction. 
\end{itemize}

\subsection{Beyond posterior consistency in Bayesian asymptotics}\label{ssec:bayesian}

\paragraph{Background:}

\iffalse
Consider the Bayesian statistical estimation framework in which given a statistical model $P_{\theta}$ a prior $\theta\sim\Pi$ is considered for the parameters of interest. 
Posterior contraction refers to the property,
$\Pi\left(L(\theta,\theta^*)>\epsilon|X\right)\rightarrow 0$
in $P_{\theta^*}$-probability. 
The main technique in the literature for proving posterior contraction with some rate is called the method of prior mass and testing \cite{barron1999consistency,ghosal2000convergence,shen2001rates,ghosal2017fundamentals}. In addition to the contraction rates, the asymptotic shape of the posterior has also been characterized by Bernstein--von Mises theorem for functional estimation in the presence of nonparametric nuisance \cite{castillo2012semiparametric,castillo2015bernstein}. 
\fi

 Under classical asymptotics, the Bernstein--Von Mises theorem guarantees that the posterior distribution concentrates around the true model at rate $1/\sqrt{n}$, and is approximately Gaussian with a center at the Maximum Likelihood Estimate (MLE). This provides a natural bridge between frequentist and Bayesian methods, and guarantees that with increasing data, the effect of the prior is washed away, and Bayesian methods match the optimal frequentist estimator. 

\iffalse
Over the past two decades, this picture has been extended to diverse high-dimensional and non-parametric models \cite{castillo2012needles,castillo2015bayesian,banerjee2021bayesian,castillo2024bayesian}. At a high-level, the picture that has emerged is as follows: if the prior is appropriately designed, the posterior will contract to the true model at the frequentist optimal rate. In addition, the Bayesian credible intervals often enjoy frequentist coverage properties, and thus allow uncertainty quantification out-of-the-box \cite{szabo2015frequentist}. 
\fi

Despite the remarkable progress in this area, most guarantees focus on the rate of contraction. Viewed through the lens of posterior contraction, many natural Bayesian methods are rate-optimal; consequently, it is difficult to compare the practical performance of competing Bayesian procedures within such a framework. This highlights the need for more refined theoretical frameworks to study Bayesian procedures.

\paragraph{Recent work:} 
A potential opportunity to obtain better theoretical understanding of Bayesian procedures is through the lens of proportional asymptotics. In the context of the linear model \eqref{eq:regression}, suppose the design entries satisfy $\bm{X}_{ij} \sim \mathcal{N}(0,1/n)$. Assume in addition that the number of features $p$ and the number of samples $n$ are both large and comparable, i.e., $p/n \to \kappa \in (0,\infty)$. 

This setup has proved to be valuable in the context of frequentist methods (e.g., see \cite{bayati2011dynamics,sur2019modern,el2013robust} and the references therein). It is useful for two reasons: (i) The framework is analytically tractable, and allows one to go beyond statistical rates and characterize the performance of statistical procedures up to constants. One may then compare the performance of competing statistical procedures via limiting behavior. (ii) Empirically, the asymptotics have excellent agreement with practical data performance. Consequently, the guarantees derived from this asymptotic framework can easily be translated into practical guidance for data analysis. 

As a concrete question, recall the proportional asymptotics setup introduced above and suppose $\bbeta_* \in [-1,1]^p$. Let $\pi$ be a prior supported on $[-1,1]^p$, and consider the corresponding posterior. What can be said about the behavior of the posterior distribution? Before proceeding further, we emphasize some differences from traditional Bayesian asymptotics. In particular, if $\bm{X}_{ij} \sim N(0,\frac{1}{n})$, standard results from random matrix theory imply that $\| \sum_{i=1}^{n} \mathbf{x}_i \mathbf{x}_i^{\top}\| = O_P(1)$ as $n \to \infty$. Consequently, the posterior does \emph{not} contract, and consistent estimation of $\bbeta_*$ is impossible. This differs sharply from traditional low-dimensional asymptotics, where $\|\sum_{i=1}^{n} \mathbf{x}_i \mathbf{x}_i^{\top}\|  = O_P(n)$. As a consequence of this scaling, the effects of the prior are not washed away even as $n \to \infty$. It might be helpful to think of this regime as a low SNR problem, where the likelihood does not overwhelm the prior effects. 

This proportional asymptotic setup has been investigated recently in information theory and statistical physics \cite{barbier2019optimal,barbier2019adaptive,lesieur2017statistical}. However, the existing literature in this context studies the \emph{matched} setting. As a concrete example, one assumes that $\bbeta_*$ is sampled from the prior $\pi$. This prior is known to the statistician, who uses the prior to construct the posterior distribution. While existing literature has focused on information-theoretic questions (e.g., the limiting mutual information between $\bbeta_*$ and the data $\mathbf{y}$), this analysis has revealed interesting properties of the posterior distribution. Under some assumptions on the prior distribution, if $\bbeta$ denotes a sample from the posterior distribution, it is understood that $\|\bbeta - \bbeta_*\|_2$ concentrates, with high probability, around a deterministic constant (dependent on the prior and $\sigma^2$). In geometric terms, the posterior mass is supported on a fixed ``band" at $O(1)$ distance from the parameter. 

\paragraph{Future directions:} From the statistical perspective, the non-matched case is more natural. In this case, one assumes that the data is sampled from a linear model \eqref{eq:regression} for some $\bbeta_* \in \mathbb{R}^p$. The statistician selects a prior $\pi$ and constructs the posterior. This leads to the natural question: what happens to the posterior distribution? One might conjecture that a similar geometric structure should emerge in this case---the posterior mass should be supported on a thin shell of radius $O(1)$. In turn, this phenomena should allow one to choose the best prior, by optimizing the distance of this shell from the model parameter $\bbeta_*$. For recent progress in the analysis of high-dimensional posteriors under mismatch, we refer to \cite{barbier2025performance,saenz2025characterizing}.

This asymptotic setup also provides a valuable framework to compare the finite-sample performance of frequentist and Bayesian procedures. The cornerstone of classical Bayesian asymptotics is the Bernstein-Von Mises theorem, which guarantees that the posterior sample is ``close" to the optimal frequentist estimator. However, the posterior sample could have better performance in finite samples---this cannot be captured through the classical asymptotic framework. In sharp contrast, under the proportional asymptotics setting introduced above, the performance of the MLE and the Bayes estimator may be computed separately and explicitly compared. This presents the intriguing possibility that an appropriate Bayes estimator might outperform the frequentist optimal estimator. Understanding the relative performance of frequentist vs.\ Bayesian methods in this high-dimensional regime is an intriguing direction for future study.

%Consider the linear regression problem \eqref{eq:regression} with $x_{ij} \sim N(0,1/n)$. Additionally, assume a proportional asymptotic regime i.e., $p/n \to \kappa \in (0, \infty)$. If the true $\bbeta_* \in [-1,1]^p$, it is natural to start with a product prior $\pi^{\otimes p}$, where $\pi$ is a prior supported on $[-1,1]$. The posterior distribution does not contract in this case. What are properties of the posterior? This is understood under the matched case in statistical physics [CITE]. Understanding this phenomenon under the non-matched case is an important problem. See recent results [CITE] for high-dimensional linear regression. In addition, this setup poses some interesting questions for the relation between high-dimensional frequentist and Bayesian inference. Specifically, is there an analogue of Bernstein-Von Mises theorem in this setting?

\subsection{Functional estimation} \label{ssec:Bayesian:functional}

\paragraph{Background:} Estimation of low-dimensional functionals is a classical question in statistical theory \cite{ibragimov2013statistical,paninski2003estimation,BickelKlaassenRitovWellner1993}. Assuming that the data arises from a family $\{\mathcal{P}_{\bm{\theta}}: \bm{\theta} \in \Theta\}$, one seeks to estimate a low dimensional functional $h(\bm{\theta}) \in \mathbb{R}$ from the observed data. For example, in the context of linear regression \eqref{eq:regression}, a natural question of interest concerns the estimation of the SNR parameter $\|\boldsymbol{\beta}_*\|_2^2$ \cite{dicker2014variance,janson2017eigenprism}. Other prominent functionals include the entropy of discrete/continuous distributions \cite{wu2016minimax}, treatment effects in causal inference \cite{robins2008higher}, estimation of linear/quadratic functionals  \cite{BickelKlaassenRitovWellner1993}, etc.  This question has been examined extensively in low-dimensional/semiparametric statistics \cite{van1991differentiable}, as well as in high-dimensions under sparsity (e.g., see the double Machine Learning literature \cite{ChernozhukovDML}). Under these assumptions, the parameter $\bm{\theta}$ can be estimated consistently, which immediately implies consistency for functional estimation. The recent progress in this direction has focused on $\sqrt{n}$ consistency, and has established gaussian fluctuations of an appropriate estimator around the truth \cite{ChernozhukovDML}.

%Recent advances have established  $\sqrt{n}$-consistent estimators with Gaussian fluctuations around the true parameter [CITE]. However, optimal estimators of low-dimensional functionals are not well-understood in high-dimensional models. 

\paragraph{Recent work:} 
Recall the linear regression model \eqref{eq:regression}. Assume the proportional asymptotics setting $p/n \to \kappa \in (0, \infty)$ and an i.i.d. Gaussian design, i.e., $\bm{X}_{ij} \sim \mathcal{N}(0, 1/n)$. As a first step, the community has focused on estimation of the regression vector $\boldsymbol{\beta}_*$. The performance of several canonical classes of estimators (e.g., penalized least squares and M-estimators) have been characterized in this setting \cite{el2013robust}. The question of optimal estimation has also been studied---specifically, the optimal estimator under certain restricted classes of estimators (e.g., generalized first-order methods or penalized convex estimators) have been characterized in the literature \cite{celentano2022fundamental,montanari2025equivalence}. 

Functional estimation is challenging in this regime, as the unknown vector $\boldsymbol{\beta}_*$ cannot be estimated consistently.
An emerging line of work focuses on functional estimation in this setup. While parameter $\boldsymbol{\beta}_*$ cannot be estimated consistently, the low-dimensional functionals can still be estimated consistently. Conceptually, the parameter $\boldsymbol{\beta}_*$ acts as a high-dimensional nuisance parameter. Functional estimation in this regime requires new ideas; we refer the interested reader to \cite{celentano2021cad,dicker2014variance} for functional estimation in the context of high-dimensional linear regression, while \cite{jiang2025new,celentano2023challenges, chen2024method} study causal effect estimation in this inconsistency regime.

  %\textcolor{red}{[Chao: could we write this problem in a Bayesian setting? It would be good to make some connection to the previous BvM content.]}

%Suppose we wish to estimate the functional $f(\bbeta_*) = \frac{1}{p} \sum_{i=1}^{p} \bbeta_{*,i}^r$. We expect that such functionals should still be estimable at the $\sqrt{n}$ rate, but the limiting variance should be larger than the classical efficient variance. Can this be characterized in general? What is a good description of the limiting variance? 

\paragraph{Future directions:}
The main question in this context concerns optimality of functional estimation in this low SNR regime. We highlight two specific examples. 
First, consider the functional $f(\bbeta_*) = \frac{1}{p} \sum_{i=1}^{p} \bbeta_{*,i}^r$.  Under classical asymptotics (i.e., $p$ fixed, $n \to \infty$), $f(\bbeta_*)$ can be estimated at a $\sqrt{n}$ rate. In addition, the optimal estimator has Gaussian fluctuations around $f(\bbeta_*)$, and the limiting variance is well-understood. However, our understanding in the proportional asymptotic setting is severely lacking. One expects that consistent and asymptotically normal (CAN) estimation should still be possible in this regime, although the limiting variance might be inflated due to the high-dimensionality of this problem. 
Another canonical example arises in the context of treatment effect estimation in causal inference. In this setting, one observes $\{(y_i, A_i, \mathbf{x}_i): 1\leq i \leq n\}$. Here, $y_i \in \mathbb{R}$ represents the observed response, $A_i \in \{0,1\}$ represents the treatment assignment, and $\mathbf{x}_i \in \mathbb{R}^p$ represents the observed covariates. Under classical assumptions in the causal literature, causal effects can be naturally related to the functional $\mathbb{E}[\mathbb{E}[y_i | A_i=1, \mathbf{x}_i]]$. Assume now a proportional asymptotic setting $p/n \to \kappa \in (0,\infty)$ and $\bm{X}_{ij} \sim \mathcal{N}(0,\frac{1}{n})$. What is the optimal estimator of this functional in this high-dimensional regime?

\iffalse
Another important challenge of Bayesian asymptotics in high-dimension is how to construct prior distributions that induce both structural assumption and natural constraint imposed from the model.

Let us revisit covariance matrix estimation in $N(0,\Sigma)$. Without any structural assumption, one would consider a conjugate prior distribution on $\Sigma$ such as inverse-Wishart. However, what if we want to estimate a sparse covariance? What would be a natural prior distribution? In sparse linear regression, one often impose the spike-and-slab prior on each regression coefficient. This prior obviously does not immediately generalize to sparse covariance. If we sample each entry of $\Sigma$ independently from a spike-and-slab prior, the joint prior of $\Sigma$ will not be supported in the class of PSD matrices. One naive modification is to condition the spike-and-slab prior on the PSD cone. However, it is unclear what the consequence is for this ad-hoc approach. Does the prior have optimal contraction rate? Does the prior lead to efficient sampling algorithm for the induced posterior?

Another example is estimation of a high-dimensional matrix that is both low-rank and sparse. Again, a natural construction of a prior distribution that respects both properties is unknown. There are many other examples such as positivity, shape-constraint and others from various problems.
\fi

\subsection{General designs and broader applicability}
\label{sec:designs_questions}

\paragraph{Background:}
Most of the work on high-dimensional statistics over the past 15 years has relied on specific modeling assumptions. In particular, in more sophisticated asymptotic regimes, such as the proportional asymptotic setting, where the number of parameters and sample size grow at the same rates, the most interpretable, sharp results have typically been established under strong structural assumptions, most notably i.i.d.\ Gaussian design matrices \cite{donoho2009message, bayati2011dynamics, thrampoulidis2018precise}. While these assumptions enable precise characterization of statistical and algorithmic phenomena, they limit the applicability of the resulting theory to more realistic data-generating processes. Consequently, relaxing such restrictive assumptions and extending high-dimensional asymptotic results to broader classes of design matrices remains one of the central open challenges in the field.

After briefly reviewing representative existing results, we will highlight several important open problems along this direction. For clarity of exposition, we focus on the canonical example of sparse linear regression introduced in Section~\ref{ssec:linear}, which serves as a setting for illustrating both the power and limitations of current high-dimensional asymptotic analyses. Suppose we observe responses $\mathbf{y} \in \R^n$ and covariates $\mX \in \R^{n \times p}$ generated from the linear model
\begin{align*}
\mathbf{y} & = \mX \bbeta_\star + \boldsymbol{\epsilon},
\end{align*}
where $\bbeta_\star$ is the unknown true coefficient vector and $\bm{\epsilon}$ represents noise. Initial results in high-dimensional asymptotics considered the case where the design matrix $\mX$ has i.i.d.\ Gaussian entries. A large body of work has studied this estimation problem under the proportional asymptotic setting where $n/p \rightarrow \delta \in (0, \infty)$ \cite{donoho2009message, bayati2011dynamics, thrampoulidis2018precise, donoho2011noise, whichbridge2020, BradicRobust2016, donoho2016high, hu2019asymptotics, bu2020algorithmic, wang2022does, su2017false}.
However, asymptotic analysis for structured, non-i.i.d.\ designs remains largely unexplored. We now highlight some of the existing challenges in this direction.

\subsubsection{Correlated, Non-Gaussian Designs} 

%A natural generalization of i.i.d.\ Gaussian designs is to remove the independence and Gaussianity assumptions. Below, we summarize some recent developments and highlight promising future directions for this area.

\paragraph{Recent work:}
Some standard tools, such as the AMP and CGMT, can be adapted to analyze dependent Gaussian designs \cite{celentano2023lasso,berthier2020state,gerbelot2023graph}. For correlated Gaussian designs, sharp risk characterizations have been derived for some popular estimators \cite{zhao2022asymptotic,celentano2023lasso,celentano2021cad,loureiro2021learning,zhang2024matrix}. In addition, several works have developed debiasing methods for constructing entry-wise confidence intervals for $\bbeta_\star$ \cite{bellec2022biasing,bellec2023debiasing}. 

Another natural generalization is a design matrix with i.i.d.\ but possibly non-Gaussian entries. A line of work has shown that under mild moment assumptions on the distribution of the entries of $\mX$, risk characterizations derived for i.i.d.\ Gaussian designs continue to remain valid, provided the mean and variance of the entries match those of the Gaussian design \cite{korada2011applications,karoui2013asymptotic,panahi2017universal,bayati2015universality, montanari2017universality, oymak2018universality, el2018impact, chen2021universality,celentano2021high,han2023universality,han2025entrywise,han2025long}. This phenomenon is called ``universality.''

More recently, these results have been generalized to the situation the samples $\bm{x}_1, \bm{x}_2, \dotsc, \bm{x}_n$ (rows of the design matrix) are independent, with possible dependence between the different features (coordinates of a particular sample $\bm{x}_i$).  In some situations, such design matrices behave like an equivalent Gaussian design with independent rows and matching row means and covariances. This is called a ``Gaussian Equivalence Principle'' in the literature \cite{hu2022universality,montanari2023universality,goldt2020modeling,mei2022generalization,gerace2024gaussian,xiao2022precise,hu2024asymptotics,ghane2024universality}. 

\paragraph{Future directions:}

  In case of the proportional asymptotics, e.g., with a correlated Gaussian model, the asymptotic risk is no longer described by a low-dimensional fixed-point equation, but rather by a high-dimensional one, which depends on the feature covariance $\bm{\Sigma}$ \cite{celentano2023lasso}. This makes it extremely challenging to compare estimators and derive useful insights from the fixed-point equations. How can we analyze the fixed-point equations and extract practically relevant information?  
    
   Another open direction is that inference methods based on the correlated Gaussian model often assume that the feature covariance matrix $\bm{\Sigma}$ is known \cite{celentano2023lasso}. Of course, in practice, $\bm{\Sigma}$ is not known and needs to be estimated. Furthermore, in high dimensions, the covariance matrix cannot be estimated consistently in the operator norm. An important problem is to develop debiasing methods that use data-driven estimators of the covariance matrix and properly account for the non-negligible estimation error in these estimators. 

In a different direction, recent exciting work in random matrix theory \cite{bandeira2023matrix} has explored a very general model of Gaussian random matrices, whose entries are jointly Gaussian, allowing for arbitrary dependence between the entries. Any such matrix $\mX$ can be written as
\begin{align*}\mX & = \sum_{i=1}^k g_i \bm{A}_i, \end{align*} for some $k \in \mathbb{N}$, where $\bm{A}_{1}, \ldots, \bm{A}_k$ are deterministic matrices and $g_{1}, \ldots, g_k$ are i.i.d.\ $\mathcal{N}(0,1)$ random variables. These works study the empirical singular value measure and  operator norm of such matrices. If $\bm{A}_1, \dotsc, \bm{A}_k$ are sufficiently non-commutative, the results lead to exact asymptotic predictions about limiting singular value measure and operator norm predictions, recovering the Marchenko-Pastur law and the Bai-Yin law in the i.i.d.\ case.  It would be very interesting to explore applications of this random matrix model in statistics. Can this model capture design matrices with interesting spatial or time series dependence? Is it possible to obtain exact risk characterizations for popular M-estimators in this model? 

Finally, understanding when such Gaussian equivalences hold, when they fail, and if they can be fixed continues to be an active research area \cite{pesce2023gaussian,xiaoyibreakdown,wen2025does}.  

\subsubsection{Rotationally Invariant Designs} 

\paragraph{Recent work:}
Under the rotationally invariant design model, the singular value decomposition of the design matrix is given by $\mX = \bm{U} \bm{\Sigma} \bm{V}^\top$, where $\bm{U}$ and $\bm{V}$ are modeled as independent, uniformly random orthogonal matrices, and the singular values $\bm{\Sigma}$ of the design are treated as deterministic. There are several reasons to study this model: First, a line of work has shown that the rotationally invariant model accurately describes the behavior of structured semi-random designs commonly seen in signal processing applications \cite{marinari1994replica,haikin2017random,donoho2009observed,monajemi2013deterministic,oymak2014case,abbara2020universality,ma2019spectral}. In statistical applications, statistical dependencies in features often result in the singular value spectrum of the design matrix deviating from the Marchenko-Pastur law. The rotationally invariant design is a cheap way to model these deviations without introducing an unknown covariance matrix, a high-dimensional nuisance parameter. Rotationally invariant designs are less well-understood than Gaussian designs. AMP algorithms for rotationally invariant designs are well-studied and are one of the only tools available to obtain exact risk characterizations for M-estimators in high dimensions \cite{ma2017orthogonal,rangan2019vector,takeuchi2019rigorous,opper2016theory,gerbelot2022asymptotic,fan2022approximate,mondelli2021pca,venkataramanan2022estimation,li2023spectrum}.  

\paragraph{Future directions:}

An interesting direction would be to develop rotationally invariant counterparts of the standard tools useful for analyzing Gaussian designs, such as:
\begin{itemize}
  \item CGMT Framework: Is there a natural generalization of various Gaussian comparison inequalities for rotationally invariant random matrices? Such a generalization could potentially lead to finite-sample versions of the currently available asymptotic risk characterizations, with simpler proofs. 
  \item Adaptive Interpolation Method: Is it possible to develop a version of the adaptive interpolation method to analyze the Bayes risk for regression problems with rotationally invariant designs? Currently, such results are only known in the low signal-to-noise ratio (SNR) regime \cite{li2023random}. 
  \item Leave-one-out method: It would also be interesting to develop a version of the leave-one-out method for rotationally invariant designs.  Such a result could be used to characterize the entry-wise distribution of debiased M-estimators and construct confidence intervals for specific coordinates of $\bbeta_\star$. The works \cite{bao2017local,lu2021householder} can provide good starting points. 
\end{itemize}

\subsubsection{Semi-Random Designs}
\paragraph{Background:} In many signal processing applications such as X-ray crystallography and Magnetic Resonance Imaging, the design matrix is highly structured with limited randomness. Two prototypical examples are as follows:
\begin{itemize}
    \item \textbf{Randomly Sub-sampled Fourier Design:} In MRI applications, the goal is to reconstruct a signal vector from its Fourier transform values at a few randomly selected frequencies. Here, the design matrix is a randomly subsampled Discrete Fourier Transform (DFT) matrix constructed by randomly selecting $n < p$ rows of a $p \times p$ DFT matrix.
    \item \textbf{Coded Diffraction Pattern (CDP) Design \cite{candes2015phase}:} In X-ray crystallography, one would like to reconstruct a signal vector from the magnitude of its Fourier transform (the phase information is lost in the measurement process). To ensure that the signal is uniquely determined, one strategy is to capture redundant measurements by imaging the same signal several times, each time modulating it with a different random mask. In this problem, the design matrix $\bm{X}$ takes the form
    \begin{align*}
        \bm{X} & = \begin{bmatrix} \bm{F}_p \bm{D}_1 \\ \bm{F}_p \bm{D}_2 \\ \vdots \\ \bm{F}_p \bm{D}_{L}  \end{bmatrix},
    \end{align*}
    where $\bm{F}_p$ denotes the $p \times p$ DFT matrix and $\bm{D}_{1}, \ldots, \bm{D}_{L}$ are random diagonal matrices representing the masks. 
\end{itemize}
Even though these matrices are quite structured, a series of works \cite{marinari1994replica,haikin2017random,donoho2009observed,monajemi2013deterministic,oymak2014case,abbara2020universality,ma2019spectral} have empirically observed that in high dimensions, the risk of many popular estimators for these structured designs is very close to the risk of these estimators for a rotationally invariant design whose singular values are identical to the singular values of the structured design. 

\paragraph{Recent work:} Proofs of this empirically observed universality result have been derived in some special cases. \cite{donoho2010counting} prove universality for the risk of the LASSO estimator for noiseless sparse linear regression. A related line of work in random matrix theory seeks to understand when  semi-random delocalized orthogonal matrices such as randomly-signed or shuffled DFT matrices behave like uniformly random (Haar) orthogonal matrices \cite{tulino2010capacity,farrell2011limiting,anderson2014asymptotically,male2011traffic,magsino2021kesten,kunisky2023generic}. More recently, \cite{dudeja2022universality,wang2024universality,dudeja2023universality,dudeja2024spectral} have characterized conditions under which the dynamics of AMP algorithms for semi-random designs (such as sub-sampled Fourier matrices or CDP matrices) are asymptotically identical to the dynamics under rotationally-invariant designs with a matching singular value spectrum.  These results yield universality principles for the risk of many popular estimators for regression problems with semi-random designs, provided that they are obtained by minimizing a \emph{strongly convex} loss function. More recently, similar results have been obtained for certain fully deterministic designs \cite{gorini2026universality}. 

\paragraph{Future directions:} In spite of  recent progress in understanding semi-random designs, many open problems remain. Existing results are restricted to strongly convex M-estimators, which precludes the popular LASSO estimator. Can we analyze convex M-estimators that are not strongly convex? Empirical simulations suggest that simple non-convex M-estimators such as spectral methods for the phase retrieval problem also exhibit universality \cite{ma2019spectral}. It would be very interesting to prove high-dimensional asymptotic risk characterizations for spectral methods for semi-random designs such as the CDP design. Lastly, we do not have any understanding of the performance of Bayes optimal estimators  for regression problems with structured designs. Do these fundamental information-theoretic limits also exhibit universality? 

%\textcolor{red}{Arian: Can we start this section more generally? Saying that semi-random matrices are important for many applications. Then, we can talk about universalities that have been observed. Then we can talk about the following mixture and Subha and Rishabh's work Zhou Fan's work? }

%\paragraph{Recent work}

  %\textcolor{magenta}{(need to double check this for accuracy)} 

%\paragraph{Sparse and Correlated Signal Recovery.} The most standard/classical model for high-dimension sparse signal recover is the Gaussian sequence model $X\sim N(\theta,I_p)$ with $\|\theta\|_0\leq k$. While almost everything is known about this model, the simple extension $X\sim N(\theta,\Sigma)$ largely remains open. Optimal estimation and testing for a sparse $\theta$ have been studied when $\Sigma$ has an equi-correlation structure \cite{kotekal2023minimax,kotekal2025sparsity}. That is, the covariance matrix has only two values depending on whether the entry is diagonal or off-diagonal. This structure corresponds to the random effect model with a rank-one component in the constant direction. For a general known $\Sigma$, such as identity plus rank-one or satisfying off-diagnoal decay, nothing is know about optimal estimation, testing and functional estimation of a sparse $\theta$. The interaction betwen sparsity and correlation is not well understood.

\subsection{Refined approximations through asymptotic expansions}\label{ssec:refinedapprox}

\paragraph{Background:}
As discussed earlier, asymptotic analysis often provides approximations for the behavior of estimators when exact finite-sample results are intractable. However, the accuracy of such approximations can depend on several factors, including the growth rate of $p$ relative to $n$  and the distributional properties of the features. Consequently, the approximations offered by the asymptotic framework may be unreliable in some settings.

A simple illustration comes from the Central Limit Theorem (CLT), one of the most fundamental tools in asymptotic analysis. For example, the sum of $n$ independent, zero-mean, unit-variance random variables, when properly normalized, converges in distribution to a Gaussian random variable with mean zero and variance one. Yet, it is well-known that this normal approximation can be inaccurate for finite sample sizes, especially when the underlying variables have heavy-tailed distributions. To address this limitation, Edgeworth proposed a higher-order expansion of the probability density function of the normalized sample mean, thereby improving the accuracy of the approximation beyond that provided by the CLT \cite{hall2013bootstrap}. Although these approximations are generally harder to obtain than the classical CLT, they yield more accurate approximations for the finite-sample distribution. 

The inaccuracy of classical asymptotic approximations has also been widely observed in the context of high-dimensional asymptotics \cite{mazumder2023subset, zheng2017does, hastie2020best}. For example, consider the canonical case of linear regression with an orthogonal design. Asymptotic minimax analysis suggests that both soft-thresholding and hard-thresholding estimators are asymptotically minimax \cite{donoho1995noising, donoho1998minimax}. However, simulation studies reveal that their finite-sample performance can differ substantially, highlighting a gap between asymptotic predictions and practical outcomes \cite{guo2023signal, hastie2020best}. These insights, when combined with inspiration from the Edgeworth expansion, prompt an important question: Can high-dimensional asymptotic analysis yield improved higher-order approximations that better capture finite-sample behavior?

\paragraph{Recent work:}
This question has not been extensively studied in the literature. However, a few partial results exist \cite{guo2023signal, ghosh2025signal, guo2024note}. We describe a simple example from one of these papers. Consider the linear regression problem~\eqref{eq:regression}, and suppose the regression coefficients belong to the following set: \begin{equation}\label{param::snr_aware}
    \Theta(k,\tau) := \Big\{ \bbeta \in \mathbb{R}^{p}: \|\bbeta\|_{0}\leq k, ~ \|\bbeta\|_{2}^{2} \leq k \tau^{2} \Big\}.
\end{equation}
Let $R(\Theta(k,\tau),\sigma)$ denote the minimax risk:
\[
R(\Theta(k,\tau),\sigma) := \inf_{\hat{\bbeta}} \; \sup_{\bbeta \in \Theta (k, \tau)} \; \mathbb{E}\bigl[\|\hat{\bbeta} - \bbeta\|_{2}^{2}\bigr]. 
\]
Note that in this model, $\frac{\tau^2}{\sigma^2}$ controls the signal-to-noise ratio. Our example is concerned with low signal-to-noise ratio. 

\begin{example}\label{thm::second-order-low-snr-minimax}
    Consider the linear regression model in Example \ref{slr:consistency:1}, with parameter space \eqref{param::snr_aware}. Suppose $k/p \rightarrow 0$ and $k/n \rightarrow 0$. Suppose further that $\tau/\sigma \rightarrow 0$. Using standard high-dimensional asymptotic analysis, one can prove that 
       \begin{equation}\label{eq:firstorderexpansion:minimax}
  R(\Theta(k,\tau),\sigma) = k \tau^{2} \big(1+o(1) \big).
    \end{equation}
However, by obtaining a higher-order expansion of the minimax risk, we can obtain  
    \begin{equation}\label{eq:secondorderexpansion:minimax}
  R(\Theta(k,\tau),\sigma) = k \tau^{2} \Big( 1 - \frac{k\tau^{2}}{p\sigma^2} \big(1+o(1)\big) \Big).
    \end{equation}
\end{example}
As is clear from equation~\eqref{eq:firstorderexpansion:minimax}, in the low-SNR regime, the first-order asymptotic approximation of the minimax risk is \(k\tau^2\), which corresponds to the performance of the zero estimator. However, in practice, many estimators can outperform the zero estimator in low-SNR settings. The second-order expansion of the minimax risk~\eqref{eq:secondorderexpansion:minimax} makes our approximation more accurate and suggests that such improvements are indeed possible. In particular, optimally tuned ridge regression can achieve both the first- and second-order terms in the minimax risk~\eqref{eq:secondorderexpansion:minimax}.

\paragraph{Future directions:}
Unlike classical asymptotics, the higher-order asymptotic approximation in high-dimensional settings is not well-understood. We list a few specific open questions here:

\begin{itemize}
\item Currently, the scope of higher-order asymptotic analysis has remained limited to a few canonical examples, such as linear regression and principal component analysis, as well as to second-order expansion results. However, we would like to develop far more general mathematical frameworks that can accommodate a broader class of statistical models and learning algorithms, thereby enabling systematic derivation of higher-order terms and deeper insights into their asymptotic behavior. 

\item A more refined asymptotic expansion of the risk of estimators can enable statisticians to differentiate between estimators that appeared equally optimal under first-order asymptotic analysis. This opens up a new avenue for studying optimal estimators through higher-order asymptotic analysis. Since higher-order asymptotic expansions provide more accurate approximations to finite-sample behavior, the conclusions drawn from such analyses about asymptotic performance can be more practically relevant.  

\item So far, higher-order asymptotic expansions have been employed to approximate the minimax risk. However, analogous questions arise naturally in other contexts, such as the Bayesian framework discussed in Section \ref{ssec:bayesian}.

\end{itemize}

%\bibliographystyle{alpha}
%\bibliography{sample}

\section{High-dimensional statistics in AI}

\subsection{Introduction}
Artificial intelligence (AI) has rapidly become one of the most transformative forces in science, technology, and society. From natural language processing and computer vision to healthcare, finance, and climate modeling, AI systems are playing central roles in decision-making and data-driven discovery. At the core of these systems lies the challenge of learning an extremely large number of parameters from vast, complex datasets. Because of the immense scale, principles of classical statistical analysis are no longer sufficient to explain or predict the behavior of modern AI models. As a result, high-dimensional statistics will play a major role in the development of AI. 

The breadth of problems and the rapid pace of development in AI make it impossible to cover all the subfields that could benefit from advances in high-dimensional statistics. Instead, we focus on a few representative subtopics that are currently of particular interest to the authors of this paper and discuss how ideas and tools from high-dimensional statistics can contribute to these areas. We hope that this selective overview will illustrate the breadth of potential interactions between the fields and encourage further cross-disciplinary collaboration.

\subsection{Fine-tuning and hyperparameter tuning}

\paragraph{Background:} Fine-tuning is the process of adapting a pretrained AI model to a specific downstream task using a smaller, task-specific dataset. Modern models are typically trained on massive, general-purpose corpora to learn broad representations, but they may perform poorly on specialized tasks or domains, especially when the target distribution differs from that of the pretraining data. Fine-tuning addresses this gap by updating some or all of the model parameters using task-specific data, allowing the model to specialize while retaining useful knowledge acquired during pretraining. Because fine-tuning datasets are usually orders of magnitude smaller than pretraining corpora, and large models may contain tens or hundreds of billions of parameters, updating all weights can be computationally and memory intensive. This has motivated the development of parameter-efficient fine-tuning methods, which adapt only a small subset of parameters while maintaining strong performance. See \cite{han2024parameter} for a review of these methods.

LoRA (Low-Rank Adaptation) \cite{hu2022lora} is the most widely adopted strategy for parsimonious fine-tuning. The key intuition behind LoRA is that, since fine-tuning data are relatively limited, the change in  weights required to specialize the model must reside within a low intrinsic-dimensional subspace. Formally, consider a pre-trained weight matrix $\bm{W}_0 \in \mathbb{R}^{d_{\text{out}}\times d_{\text{in}}}$. In standard fine-tuning, $\bm{W_0}$ is adjusted to $\bm{W}' = \bm{W}_0 + \bm{\Delta W}$, where $\bm{\Delta W}$ is the update matrix computed via backpropagation. LoRA constrains $\bm{\Delta W}$ by reparametrizing it as a product of two low-rank matrices:
\[
\bm{W}' = \bm{W_0} + \gamma \bm{B} \bm{A},
\]
where $\bm{B} \in \mathbb{R}^{d_{\text{out}}\times r}$ and $\bm{A} \in \mathbb{R}^{r\times d_{\text{in}}}$, and $\gamma$ is a constant scaling factor (often set relative to $r$). The hyperparameter $r$, known as the rank, is chosen such that $r \ll \min(d_{\text{out}}, d_{\text{in}})$. By construction, the update matrix $\bm{\Delta W} = \gamma \bm{B} \bm{A}$ satisfies $\text{rank}(\bm{\Delta W}) \leq r$. In effect, LoRA imposes a structural bottleneck in the forward pass so that the influence of any input on the weight update must pass through an $r$-dimensional latent space. This low-rank structure is different from sparsity-based regularization (e.g., $\ell_1$ or LASSO), which assumes many individual parameters are exactly zero. In contrast, the LoRA update $\bm{\Delta W}$ is generally a dense matrix, albeit one that is structurally constrained by its factorization.

During fine-tuning, LoRA keeps the original weights $\bm{W_0}$ frozen and updates only the matrices $\bm{A}$ and $\bm{B}$.
 The reduction in trainable parameters is substantial: instead of optimizing $d_{\text{out}} \cdot d_{\text{in}}$ parameters, LoRA optimizes only $r(d_{\text{out}} + d_{\text{in}})$ parameters. To ensure the adaptation begins smoothly, $\bm{A}$ is typically initialized with random Gaussian entries, while $\bm{B}$ is initialized to zero. At inference time, the merged matrix $\bm{W}'$ is used, resulting in no additional computational latency compared to the original model. While LoRA drastically reduces the memory footprint and storage costs by design, it has been empirically demonstrated to achieve performance on par with full fine-tuning.

\paragraph{Recent work:} From a statistical perspective, LoRA offers a compelling framework for structured dimension reduction in the parameter space, aligning closely with concepts prevalent in high-dimensional statistics. LoRA can be viewed as performing an implicit PCA on the hypothesized full update matrix $\bm{\Delta W}$, assuming that the necessary changes to the weights are concentrated along a small number of principal directions in the parameter space. Unlike traditional PCA applied to data features, however, LoRA applies this low-rank constraint to the \textit{change} in the parameter matrix and learns this structure directly via optimization on the fine-tuning objective. 

%\textcolor{red}{AM: there are a few papers that have some theoretical results here: "LoRA Training Provably Converges to a Low-Rank Global Minimum or It Fails Loudly (But it Probably Won’t Fail)", "Understanding the Learning Dynamics of LoRA: A Gradient Flow Perspective on Low-Rank Adaptation in Matrix Factorization", "LoRA-One: One-Step Full Gradient Could Suffice for Fine-Tuning Large Language Models, Provably and Efficiently" "Why LoRA Fails to Forget: Regularized Low-Rank Adaptation Against Backdoors in Language Models". Should we cite these papers? Can we mention something about a few other algorithms such as MeZO?}

Beyond this heuristic interpretation, a rapidly growing body of theoretical literature has begun to analyze optimization and statistical properties of LoRA. From an optimization viewpoint, \cite{kim2025lora} established that LoRA training either converges to a low-rank global minimum or becomes trapped in a high-rank, large-norm bad minimum. Small standard LoRA initialization is found to produce a two-phase dynamic that reaches an error floor set by singular-space misalignment \cite{xu2025understanding}. Additionally, LoRA's updates are shown to align with the top singular subspaces of the one-step full fine-tuning gradient \cite{zhang2025lora}. On the negative side, \cite{luong2026lora} demonstrated that LoRA fails to remove backdoors, mainly because its low-rank updates are too weak and spectrally misaligned to override poisoned pretrained directions.

\paragraph{Future work:}
The connections between LoRA and statistical principles suggest many opportunities for statisticians, particularly those specializing in high-dimensional inference, to contribute to the theoretical understanding and methodological improvement of parameter-efficient fine-tuning methods:

\begin{enumerate}
   \item \textbf{Theoretical Characterization of Hyperparameter Tuning.} The optimal configuration for LoRA remains unclear and is largely determined through trial and error. For instance, extensive experiments by Schulman and colleagues at Thinking Machines Lab have shown that LoRA often requires a significantly higher learning rate than full fine-tuning, and applying LoRA broadly across different layer types (e.g., both attention and MLP modules) is generally necessary to match the capacity of full fine-tuning. Developing a rigorous theory that elucidates, e.g., the relationship between LoRA and full fine-tuning learning rates would be highly valuable for guiding practice.

    \item \textbf{Adaptive Rank Selection (Model Selection).} Selecting the rank $r$ involves a fundamental trade-off between computational efficiency and model fidelity. The optimal $r$ is governed by the \textit{intrinsic dimension} of adaptation, which is influenced by task complexity, dataset size, and the divergence between pre-training and fine-tuning distributions. Formalizing these relationships remains an open challenge. Statisticians can contribute by developing adaptive rank selection methods, leveraging tools from random matrix theory, Bayesian nonparametrics, or tailored information criteria (e.g., AIC/BIC analogs) to estimate this intrinsic dimensionality.

    \item \textbf{Structured Low-Rank Methods.} Researchers could explore structured LoRA variants. For instance, imposing sparsity or group-sparsity constraints on the factors $\bm{A}$ and $\bm{B}$ could lead to even more efficient,  interpretable adaptations. One could also consider adding regularization on certain norms of the matrices $\bm{A}$ and $\bm{B}$, thereby combining the benefits of low-rank and sparse structures. This approach is reminiscent of the extensive literature in high-dimensional statistics on structured sparsity, and techniques and insights from that domain may prove valuable for advancing LoRA methodology.

\item \textbf{Training Dynamics and Optimization Landscape.} Theoretical understanding of LoRA's training dynamics remains limited. Analyzing the non-convex optimization landscape induced by the $\bm{B}\bm{A}$ factorization is crucial, particularly for elucidating how the low-rank structure affects convergence properties, implicit regularization, and ultimately generalization performance. Developing a rigorous theory that connects the geometric properties of this factorized optimization landscape to LoRA's empirical success would be valuable for optimizing training regimes and deepening our understanding of the method's inductive biases.

\end{enumerate}

\subsection{In-context learning}\label{ssec:incontext}

\paragraph{Background:} Large language models have enabled a fundamentally new learning paradigm that differs markedly from traditional supervised fine-tuning. In \emph{in-context learning} (ICL), a frozen pre-trained model performs new tasks by conditioning its output on a prompt consisting of a few input-output examples, without any gradient updates or model parameter changes \cite{mccann2018natural, radford2019language, brown2020language}. For instance, given a prompt containing several (English phrase, French translation) pairs, the model can translate a novel English phrase—despite never having been explicitly trained for translation as a downstream task.

This phenomenon is surprising for several reasons. First, LLMs are trained solely on next-token prediction, yet at test time, they exhibit competence on a wide range of downstream tasks that were never explicitly part of the training objective. Second, the model weights remain entirely fixed during ICL; the ``learning" occurs purely through the structure of the prompt. Third, ICL was not apparent in smaller models, but appears to \emph{emerge} as models scale, raising fundamental questions about the relationship between model capacity, training data diversity, and generalization.

From a practical standpoint, ICL enables rapid prototyping of new tasks without large-scale data collection, expensive compute for fine-tuning, or specialized know-how. From a scientific standpoint, it represents a new, exciting, and poorly understood learning paradigm that challenges conventional statistical thinking about generalization.

\paragraph{Recent work:}
To complement the extensive empirical literature on ICL, which often yields contradictory conclusions due to sensitivity to experimental settings and hyperparameters, the theoretical ML community has sought tractable frameworks for rigorous analysis. A particularly influential framework, introduced by Garg et al.\ \cite{garg2022can}, situates ICL within the classical estimation setting familiar to the rest of this paper. Consider a collection of $k$ linear regression tasks, each parametrized by a $d$-dimensional coefficient vector $\bbeta_j \in \mathbb{R}^p$, for $j = 1, \ldots, k$. For each task, we generate a training sequence of $T+1$ input-output pairs:
\[
\left\{(\bx_1, y_1), (\bx_2, y_2), \ldots, (\bx_T, y_T), (\bx_{T+1}, y_{T+1})\right\},
\]
where $y_t = \bbeta_j^\top \bx_t + z_t$ and $z_t$ represents noise. The model is trained on $n$ such sequences, with the objective of predicting $y_{T+1}$ given the \emph{prompt} $(\bx_1, y_1, \ldots, \bx_T, y_T, \bx_{T+1})$.

Within this framework, ICL manifests in two distinct modes \cite{pan2023context}. \emph{Retrieval-based ICL} corresponds to in-distribution generalization: the model correctly handles sequences generated from one of the $k$ task vectors seen during pre-training, even if the specific sequence was not encountered. In contrast, \emph{true ICL} corresponds to out-of-distribution generalization: the model successfully learns from a sequence generated by a novel task vector that has not been among the $k$ pre-training tasks. Within the linear regression ICL setup described above, Raventos et al.~\cite{raventos2023pretraining} empirically demonstrated that transformers can perform both modes of ICL, and remarkably, their estimation error matches that of the least-squares estimator—suggesting that pre-trained transformers behave as Bayes-optimal predictors in this setting \cite{xie2022explanation}.

These observations have catalyzed a rich body of theoretical work on the linear regression 
ICL setting and its natural generalizations, including ICL over tasks defined by generalized 
linear models. One line of research seeks to understand how different architectures can perform ICL through explicit construction of weights that implement the task \cite{von2023transformers, akyurek2022learning,edelman2024evolution,li2024fine,li2023transformers}. A still-debated hypothesis posits that each transformer layer implements one step of gradient descent, so that with sufficient depth the model approximates the least-squares solution \cite{von2023transformers,shen2024position,bai2023transformers,ahn2023transformers,mahdavi2024revisiting,fu2024transformers}.  Other work examines how the diversity of the pre-training distribution affects ICL capabilities, how far test-time tasks can deviate from pre-training tasks while still enabling successful ICL, and the precise role of model scale relative to task difficulty \cite{raventos2023pretraining,park2024competition,carroll2025dynamics,deora2025context}.

From a statistical perspective, linear regression ICL introduces a compelling 
twist on classical supervised learning that is naturally suited to high-dimensional analysis.  The data structure involves four key parameters: the number of training sequences $n$, the input dimension $p$, the number of pre-training tasks $k$, and the context length $T$. A model mapping sequences of $p$-dimensional vectors and their labels to predictions requires $O(p^2)$ learnable parameters in the simplest (linear attention) case. The natural high-dimensional asymptotic regime is thus $n \asymp p^2$, $T \asymp p$, and $k \asymp p$. Recent work \cite{lu2025asymptotic,letey2025pretrain} has derived sharp asymptotic predictions in this proportional regime for a simplified single-layer linear attention model. In essense, the prediction takes the form $\langle \bx_{T+1}, \bm{\Theta} \hat{\bbeta}_\text{avg} \rangle$ with $\hat{\bbeta}_\text{avg} = \sum_{t=1}^{T} \bx_t y_t$ (the averaging estimator from the prompt) and $\bm{\Theta}$ the learnable parameter matrix. %\textcolor{red}{I am not sure if I understand the notation in the last sentence. What do we do with $w$? }

\paragraph{Future directions:}
High-dimensional statistics has an important role to play in addressing the many open questions surrounding ICL:

\begin{enumerate}[leftmargin=*]
    \item \textbf{Sharp asymptotic and finite-sample analysis of ICL generalization.}  Beyond the two recent works mentioned above, deriving precise asymptotic formulas for both in-distribution (retrieval) and out-of-distribution (``true'' ICL) generalization error in the proportional regime remains largely open, even for simplified architectures. Understanding how these errors depend on $n$, $p$, $k$, $T$, and properties of the input data would provide  insight into ICL mechanisms and could provide guidance for large-scale empirical studies.

    \item \textbf{Tradeoffs Between ID and OOD ICL.} There appears to be a fundamental tension between retrieval-based and true in-context learning, e.g., \cite{pan2023context,lin2024dual,raventos2023pretraining,lu2025asymptotic}. How do these tradeoffs depend on model size, context size, pre-training duration, and task diversity? Characterizing these tradeoffs rigorously would clarify when and why ICL succeeds or fails.

    \item \textbf{Algorithm and pre-training duration dependence.}
Most existing theories either assume an idealized Bayes predictor or construct pre-training 
weights capable of performing ICL, without guarantees that standard training procedures 
actually recover such weights. Exceptions require strong simplifying assumptions \cite{zhang2024trained}. In practice, ICL emerges on transformer-like architectures, under finite compute budgets and 
next-token-prediction pre-training objectives. Developing algorithm-dependent guarantees 
that characterize ICL performance as a function of optimization time, algorithm, and 
architecture would bridge this gap between theory and practice.

    % \item \textbf{Role of Regularization in Pre-training.} Regularization techniques (explicit or implicit through choice of optimizer) during pre-training likely influence ICL capabilities, yet their effects are poorly understood. Developing theory that connects pre-training regularization to downstream ICL performance is an important open direction.

    \item \textbf{Beyond Linear Models.} Most rigorous theoretical results concern linear 
regression tasks and/or single-layer linear attention architectures. Extending the analysis 
to structured regression (e.g., sparse or low-rank coefficients), classification tasks, 
mixture models, and, even more importantly, discrete sequential data such as Markov chains would significantly 
broaden the applicability of the theory and help close the gap between existing theory 
and practice.

    \item \textbf{Effect of Scale Relative to Task Complexity.} The emergence of ICL in larger models suggests a phase transition phenomenon. Precisely characterizing the scaling laws (cf.\ Section \ref{sec:scaling_laws}), i.e., how model size, data size, and task complexity interact to enable or preclude ICL, is a natural direction for high-dimensional analysis of ICL.
\end{enumerate}

\subsection{Machine unlearning}

\paragraph{Background:} Recent regulations, such as the California Consumer Privacy Act (CCPA) \cite{pardau2018california} and the General Data Protection Regulation (GDPR) \cite{mantelero2013eu}, grant individuals the ``right to be forgotten," allowing them to request the deletion of personal data from various organizational databases. While removing personal data from a training dataset is straightforward, the real challenge arises when these data points have already been used in the training of models. Companies are therefore required not only to remove an individual datapoint, but also eliminate any trace of the data from all relevant models, which can be very costly. In response to this challenge, the area of machine unlearning \cite{cao2015towards, bourtoule2021machine, nguyen2025survey} has emerged. More formally, consider a dataset
\[
\mathcal{D} = \left\{(\bx_1, y_1), (\bx_2, y_2), \ldots, (\bx_n, y_n)\right\},
\]
where each feature vector $\bx_i \in \mathbb{R}^p$ has corresponding response $y_i \in \mathbb{R}$. 
The goal is to predict $y_i$ based on $\bx_i$. 
Suppose we use a parametric model $g_{\bm{\theta}}(\bm{x})$, parametrized by $\bm{\theta} \in \mathbb{R}^d$, to map input features to responses. 
An estimate of $\bm{\theta}$ is obtained by minimizing the \emph{empirical risk}:
\[
\bm{\hat{\theta}} \in \arg\min_{\bm{\theta}} \sum_{i=1}^n \ell(y_i, g_{\bm{\theta}}(\bx_i)),
\]
where $\ell(\cdot, \cdot)$ denotes the loss function, e.g., the squared loss for regression or the cross-entropy loss for classification.
 Suppose a request has been made to remove a subset of data points, denoted by 
$\mathcal{M} \subset \mathcal{D}$, from the dataset. 
The goal is to recompute the model parameters after excluding these data points:
\begin{equation}\label{eq:dataremoval:opt}
\bm{\hat{\theta}_{\setminus \mathcal{M}}} \in 
\arg\min_{\bm{\theta}} \sum_{i \notin \mathcal{I}_{\mathcal{M}}} 
\ell(y_i, g_{\bm{\theta}}(\bx_i)),
\end{equation}
where $\mathcal{I}_{\mathcal{M}}$ denotes the indices of the data points 
belonging to $\mathcal{M}$.

The most straightforward approach is to solve the above optimization problem 
from a random initialization. However, this can be computationally expensive, 
especially for large models and datasets. A more efficient alternative is to solve 
\eqref{eq:dataremoval:opt} using an iterative optimization algorithm 
initialized at the previously trained parameter estimate $\bm{\hat{\theta}}$. 
Since the removal set is typically small (i.e., $|\mathcal{M}| \ll |\mathcal{D}|$), 
the new optimum $\bm{\hat{\theta}_{\setminus \mathcal{M}}}$ is expected to be close to 
$\bm{\hat{\theta}}$, allowing convergence within only a few iterations.

While this approach is computationally efficient, it carries the risk that the 
resulting model may still retain residual information about the removed 
data points, leading to privacy leakage. To mitigate this, methods 
inspired by privacy-preserving learning literature introduce a small 
amount of random noise during the update process to obscure any remaining 
information related to $\mathcal{M}$. Suppose  the outcome of the machine unlearning algorithm is $\bm{\tilde{\theta}_{\setminus \mathcal{M}}}$. Two requirements are expected from this approximation \cite{guo2019certified, dwork2006differential, zou2025certified, pandey2025gaussian}:

\begin{enumerate}
    \item (Certifiability.) No information about the data points in $\mathcal{M}$ should be recoverable from the unlearned model, as the users have requested that their data be completely removed.

    \item (Accuracy.) $\bm{\tilde{\theta}_{\setminus \mathcal{M}}}$ is expected to approximate $\bm{\hat{\theta}_{\setminus \mathcal{M}}}$ in downstream metrics such as out-of-sample prediction accuracy.

\end{enumerate}

\paragraph{Recent work:}
The work in the field of machine unlearning falls into the following categories:
\begin{itemize}
\item Algorithmic approaches: A range of algorithms have been proposed for machine unlearning under a wide range of settings. For instance, different authors have considered the application of a Newton method and its approximations (e.g., based on an approximation of the inverse Hessian) and gradient-based methods \cite{guo2019certified, zou2025certified, sekhari2021remember, neel2021descent, pandey2025gaussian}. In addition to these methods, several ideas have been developed to make the computations faster, such as SISA \cite{bourtoule2021machine}and TRAK \cite{park2023trak}. A recent survey of the existing machine unlearning algorithms can be found in \cite{nguyen2025survey}. 

\item Several papers have studied the accuracy and certifiability of certain machine unlearning models \cite{guo2019certified, zou2025certified, sekhari2021remember, neel2021descent, pandey2025gaussian}. However, the assumptions adopted in most of these works hold only in low-dimensional settings. See \cite{zou2025certified, pandey2025gaussian} for examples of major discrepancies that can arise between low-dimensional and high-dimensional analyses of machine unlearning algorithms. 
\end{itemize}

\paragraph{Future directions:}

As mentioned above, most theoretical work in the machine unlearning literature has focused on low-dimensional settings and convex training objectives. However, many modern models operate in high-dimensional, often nonconvex regimes, creating a need to understand the accuracy and limitations of unlearning methods in these settings. A few recent papers have begun to extend the analysis to high-dimensional or nonconvex scenarios \cite{zou2025certified, tong2026imperfect, pandey2025gaussian}, typically by studying specific algorithms and their performance. Nevertheless, the accuracy and certifiability of most unlearning methods remain poorly understood, and the fundamental limits of machine unlearning are  largely unknown.  

 Another central challenge in machine unlearning---and more broadly in  privacy---is the empirical verification that publicly released models or estimates do not reveal information about users who have requested the removal of their data. A growing body of work has explored how to empirically  validate privacy guarantees. Early studies examined membership inference and hypothesis testing attacks as practical tools for detecting privacy leakage and evaluating whether DP mechanisms meaningfully limit adversarial advantage \cite{dwork2015reusable, yeom2018privacy}; subsequent research developed empirical frameworks for auditing the privacy of machine learning algorithms \cite{jayaraman2019evaluating, mokander2024auditing, steinke2023privacy}. However, a major gap remains: the reliability of these auditing methods in high-dimensional settings is not well understood, nor are the most efficient strategies for auditing machine unlearning algorithms.

% \begin{itemize}
% \item 

% \item 

% \item 

% \item 

% \item 

% \end{itemize}

%reduced-rank regression

\subsection{Mechanistic interpretability}

\paragraph{Background:}
The black-box nature of LLMs, when deployed in high-stakes applications, has led to growing concerns around interpretability, safety, alignment, and trustworthiness. While the pursuit of interpretable machine learning has a rich history---ranging from the stability-driven framework \cite{yu2020veridical} to  network dissection work \cite{bau2020understanding}---the  subfield of mechanistic interpretability has recently gained significant momentum. 

Popularized largely by researchers at OpenAI and Anthropic (e.g., \cite{nanda2023mechanistic}), this approach seeks to reverse-engineer the exact computations performed inside neural networks \cite{olah2020zoom}. Rather than treating models as black boxes and analyzing input-output correlations, mechanistic interpretability aims to construct an algorithmic understanding of how components (e.g., attention heads, neurons, weight matrices) interact. The ultimate goal is to translate trained model weights into interpretable pseudocode or circuit diagrams. In statistical terms, viewing an LLM as a black-box function with billions of parameters, mechanistic interpretability attempts to decompose the function into a network of subfunctions and latent variables, akin to analyzing a complex high-dimensional estimator by identifying meaningful intermediate variables and their relationships.

For a simple example, consider a transformer with $L$ layers, where each layer $\ell$ contains $H$ attention heads and an MLP module. The model processes an input sequence of tokens $\mathbf{x} = (x_1, \ldots, x_n)$ by iteratively updating a sequence of hidden representations $\mathbf{h}^{(\ell)} = (\bm{h}_1^{(\ell)}, \ldots, \bm{h}_n^{(\ell)})$, where $\bm{h}_i^{(\ell)} \in \mathbb{R}^{d}$ denotes the representation of token $i$ at layer $\ell$ (often called the residual stream) \cite{elhage2021mathematical}. At each layer, attention heads compute weighted combinations of previous token representations based on learned query-key-value projections, while MLP layers apply position-wise nonlinear transformations. Mechanistic interpretability seeks to identify which specific heads or neurons are responsible for particular behaviors. For instance, a circuit for pronoun resolution might consist of a subset of attention heads $\{h_{\ell_1, j_1}, h_{\ell_2, j_2}, \ldots\}$ across various layers that attend to antecedent candidates, combined with specific MLP neurons that integrate this information to boost the probability of the correct referent \cite{wang2022interpretability, olsson2022incontext}. Formally, if we denote the output of head $h_{\ell,j}$ as $\mathbf{o}_{\ell,j} \in \mathbb{R}^{n \times d}$ and the MLP output at layer $\ell$ as $\mathbf{m}_{\ell} \in \mathbb{R}^{n \times d}$, a circuit hypothesis posits that a particular behavior (e.g., correctly predicting a pronoun's referent) can be attributed to a sparse subset $\mathcal{S} \subset \{(\ell, j) : 1 \leq \ell \leq L, 1 \leq j \leq H\} \cup \{(\ell, \text{MLP}) : 1 \leq \ell \leq L\}$ of components.

\paragraph{Recent work:} At a high level, mechanistic interpretability bears substantial similarity to high-dimensional statistics in that both seek to select interpretable functionalities (circuits or variables) from an ocean of many possible candidates. Indeed, researchers have developed techniques such as probing classifiers \cite{alain2016understanding}, which serve as diagnostic tests and connect to classical estimation theory: training a linear probe on activations is analogous to fitting a regression model to test whether a particular feature is linearly separable in the representation space, with probe performance providing evidence for or against the null hypothesis of independence between activations and the target feature. Another technique, feature attribution methods, quantifies the marginal contribution of individual neurons or heads to outputs, thereby guiding the search for important components. Linear representation analysis \cite{burns2023discovering,uppaal2024model,veitch2024linear} investigates whether high-level concepts correspond to interpretable linear directions in activation space, often employing dimensionality reduction techniques such as PCA or sparse dictionary learning to disentangle polysemantic neurons (neurons that respond to multiple unrelated features due to superposition) \cite{elhage2022toy,song2023uncovering}. Recent advances include the use of sparse autoencoders \cite{bricken2023towards} to factorize high-dimensional activations into overcomplete bases where individual features are more semantically atomic, yielding tens of millions of candidate interpretable features in analyses.

\paragraph{Future directions:}
These deep connections between mechanistic interpretability and statistical methodology suggest several promising research directions, and opportunities are abundant. A challenge is the lack of objective, quantitative measures of interpretability quality \cite{sharkey2025open}. While circuit hypotheses can be validated through interventions \cite{goldowsky2023localizing}, there is no consensus on how to assess whether an interpretation is complete, minimal, or generalizable. Statisticians could contribute by developing principled evaluation frameworks—for instance, defining notions of ``variance explained'' by a proposed circuit (analogous to $R^2$ in regression), or formulating information-theoretic measures that quantify how much of the model's decision-making is captured by identified components. Another direction is to develop scalable, data-driven methods for automatically discovering circuits \cite{sharkey2025open}, which would dramatically accelerate progress. Given observational and interventional data (activations under various inputs and ablations), the problem is to infer the graph of interactions among components. Techniques from high-dimensional graphical modeling, sparse regression with structured penalties (e.g., group Lasso for identifying clusters of related neurons), or Bayesian structure learning could be adapted to this setting. Moreover, as the field matures, ensuring that interpretable discoveries are reproducible and statistically sound
becomes critical. Many current findings are based on small sets of examples or qualitative assessments, raising concerns about cherry-picking or overfitting interpretations to specific inputs \cite{rauker2023toward}. Statisticians can develop multiple hypothesis testing frameworks for circuit validation when screening thousands of neurons.

\subsection{Neural scaling laws}\label{sec:scaling_laws}

\paragraph{Background:}
When training modern deep learning models, three key quantities govern the efficacy of training: the \emph{model size} $p$ (i.e., the number of trainable parameters), the \emph{dataset size} $n$ (i.e., the number of training examples), and the \emph{number of training iterations} $S$ (i.e., the total number of gradient updates applied to the model weights). The relationship between $S$ and $n$ is controlled by the \emph{number of epochs} $E$ (how many times the model sees each example during training) and the \emph{batch size} $B$ (how many examples are used to form each gradient update), satisfying $nE = BS$. We denote by $m = BS = nE$ the total number of examples---potentially with repetitions---visited during training. A basic (and widely used) abstraction is that training compute, measured in FLOPs, scales proportionally to the product of model size and the total number of processed examples
$
\mathrm{FLOPs} \ \propto \ p\, m,
$
up to architecture-dependent constants \cite{kaplan2020scaling}.
This framing leads to a central design question: \emph{for a fixed compute budget, how should one allocate resources across $p$ and $m$ (equivalently, $n$ and $E$, or $B$ and $S$) to optimize generalization performance?}

In the early era of LLM pretraining, compute constraints typically precluded multiple passes over the data, making the single-epoch regime ($E = 1$) a natural setting. With batch size $B$ fixed, the question reduces to determining the optimal trade-off between $p$ and $n$ to minimize population risk---measured as next-token cross-entropy loss on a held-out test corpus (often reported as \emph{perplexity} in the NLP literature)---subject to a compute constraint, which in this regime amounts to controlling $p\cdot m = p\cdot n$. Notably, since pretraining uses a self-supervised objective (next-token prediction), increasing $n$ corresponds to using a larger text corpus without explicit labeling.

\paragraph{Recent work:}
Empirical studies have found remarkably regular \emph{neural scaling laws} for language modeling loss as a function of model  and dataset size \cite{kaplan2020scaling, hoffmann2022training}. Over a wide range of large $p$ and $n$, the test loss can often be fit by a power-law form such as
\[
\mathcal{L}(p,n) \approx H + c_p\, p^{-a_p} + c_n\, n^{-a_n},
\]
where $H$ is an irreducible term (often interpreted as an entropy-like floor for the data distribution), and the constants $(c_p,c_n,a_p,a_n)$ are obtained by fitting learning curves. These fitted laws provide a practical mechanism for compute-efficient scaling: given a target risk, they suggest how to trade off model size and data, and conversely, given a compute budget, they suggest approximately compute-optimal training configurations \cite{kaplan2020scaling, hoffmann2022training}. This methodology has influenced the development of many major LLMs, both commercial and open-source \cite{touvron2023llama, brown2020language, jiang2023mistral}. Importantly, the fitted constants are not universal: they depend on architecture, data distribution and curation (e.g., tokenization), and the optimization procedure (including learning-rate schedules and other hyperparameters).

\paragraph{Future directions:}
Neural scaling laws raise profound theoretical questions. The experiments required to calibrate scaling laws are extraordinarily costly, limiting the ability of the broader research community to conduct systematic empirical studies of the origins, mechanisms, limitations, and dependencies of scaling behavior. \emph{What properties of the data, task, or architecture give rise to power-law scaling? Under what conditions do scaling laws break down? How sensitive are the exponents $a_p$ and $a_n$ to distributional properties of the training data or architectural choices?} These questions remain largely open.

Moreover, the regime of practical relevance is shifting. As computational resources have grown, the AI community has approached a regime where much of the available high-quality data has been exhausted: training corpora are already at the scale of the accessible internet \cite{villalobos2022will,muennighoff2023scaling}. We are thus transitioning from a \emph{compute-constrained} regime to a \emph{data-constrained} regime, where the single-epoch assumption ($E = 1$) may no longer be appropriate. New questions emerge: \emph{How many epochs should one use? How does repeated exposure to data affect generalization, and at what point do diminishing returns or overfitting set in?} Other optimizer hyperparameters (e.g., batch size, learning-rate schedules, momentum coefficients, weight-decay, regularization) are also becoming increasingly important as model performance approaches saturation and marginal gains require finer tuning.

This is precisely where statistical theory can play an informative role: not only explaining the mechanisms underlying existing empirical observations, but also pushing the boundaries of analysis into regimes that may not yet be practically relevant but will become so as compute and data landscapes evolve. The fact that scaling laws are most clearly observed at large $p$ and $n$ suggests that proportional high-dimensional regimes are natural mathematical limits in which to seek principled predictions. Indeed, a recent line of theoretical work has begun to investigate scaling laws mathematically \cite{bahri2024explaining,maloney2022solvable,paquette20244+,lin2024scaling,yan2025larger,lin2025improved,renemergence,arouslearning}, but we are only at the surface, and the pace of empirical progress demands accelerated theoretical development. Several concrete directions merit attention:

\begin{enumerate}

\item \textbf{Beyond simple models.} The rigorous theory of scaling laws has progressed from linear regression to random feature models, and very recently to models that incorporate feature learning in shallow networks. Substantial gaps remain as, for example, none of these analyses address the transformer architecture or attention mechanisms---the very models for which empirical scaling laws were originally documented. More broadly, extending the theory to deep architectures and richer task structures (e.g.o, classification, structured prediction, sequential models) that better reflect language modeling are important open directions. A related puzzle is explaining the empirical observation that fitted exponents appear to depend more on properties of the data distribution than on specific architectural details, a ``universality'' phenomenon only partially accounted for by current theory.

    \item \textbf{The data-constrained regime.} With the transition from compute-constrained to data-constrained training, there is an urgent need for theory that addresses multi-epoch training. How does generalization error evolve as a function of the number of epochs $E$? What is the optimal allocation between $n$ and $E$ for fixed $m=nE$? How do memorization, curriculum effects, and data ordering interact with scaling laws?

    \item \textbf{Algorithm dependence.} Scaling behavior depends on optimization details (e.g., batch size, schedules, regularization, data ordering). Developing theory that predicts how fitted exponents and constants vary with algorithmic choices is an important gap. A related question is whether there exist regimes in which interpolation yields \emph{benign overfitting} phenomena in language modeling, analogous to those in supervised learning \cite{bartlett2020benign,belkin2019reconciling}.

    \item \textbf{Scaling laws for emerging paradigms.} The original scaling laws were formulated for standard pretraining via next-token prediction. New training and inference paradigms are rapidly emerging. What are the scaling laws for \emph{in-context learning} (Section~\ref{ssec:incontext}) as a function of model size, context length, and number of in-context examples? Similarly, \emph{test-time compute} (e.g., additional computation at inference via sampling or multi-step reasoning) is becoming an important lever for performance \cite{snell2024scaling}; characterizing scaling laws that jointly account for training compute and test-time compute is an increasingly important direction.
\end{enumerate}

\subsection{Reinforcement Learning with Verifiable Rewards}

As pre-training data has become increasingly saturated, much of the focus in LLM development 
has shifted from pre-training to post-training. A simple and empirically powerful form of 
post-training is \emph{reinforcement learning with verifiable rewards} (RLVR), especially in 
domains such as mathematics and coding, where correctness can be checked efficiently by an 
external oracle \cite{RLVR_term,r1}. Formally, given a distribution $\mathcal{P}$ over 
(problem, answer) pairs $(x,a)$, and a model distribution $\pi_{\bm{\theta}}(\cdot\mid x)$, 
RLVR updates the high-dimensional parameter vector $\bm{\theta}$ to maximize
\[
\mathbb{E}_{(x,a)\sim\mathcal{P}}\bigl[\rho_{\bm{\theta}}(x,a)\bigr], \qquad
\rho_{\bm{\theta}}(x,a)=\mathbb{E}_{y\sim\pi_{\bm{\theta}}(\cdot\mid x)} r(y,a),
\]
where $r(y,a)$ is typically a binary reward indicating correctness. In practice, RLVR uses 
policy-gradient updates based on sampled responses and corresponding verified rewards, with 
algorithms such as REINFORCE \cite{williams1992simple}, RLOO \cite{kool2019buy}, and GRPO \cite{GRPO_original} differing mainly in how they define the advantage 
weights.

We see an opportunity for high-dimensional statistics to play a central role in the 
theoretical understanding of RLVR, precisely because the field has developed tools for 
the regimes that make RLVR challenging: highly overparameterized models, sparse and 
compressed feedback, finite-sample effects in gradient estimation, and training samples 
whose distribution depends on the current iterate. Concretely, the parameter dimension is large relative to the amount of  supervision, the feedback is often binary and outcome-based, the policy-gradient signal must be estimated from a small number of sampled responses per prompt, and these responses are 
generated by the current policy, so their distribution changes throughout learning.

\paragraph{Recent work:}
Recently, \cite{thrampoulidis2025advantage,davis2025objective} has shown that, for binary rewards, many RLVR 
algorithms can be interpreted as stochastic gradient ascent on a surrogate-reward objective
\[
\mathbb{E}_{(x,a)\sim\mathcal{P}}\bigl[F(\rho_{\bm{\theta}}(x,a))\bigr]
\]
for a suitable surrogate $F$. Under this perspective, REINFORCE and RLOO correspond to a 
linear surrogate, while GRPO corresponds to $F(\rho)=2\arcsin\sqrt{\rho}$. This provides a 
natural bridge to high-dimensional statistics: different choices of $F$ define different 
estimators under sparse binary feedback, and therefore induce different bias-variance tradeoffs, 
different emphasis on easy versus hard prompts, and potentially different statistical 
behavior. Supporting this perspective, \cite{heckel2026asymmetric} showed that the log-surrogate 
$F(\rho)=\log\rho$ can outperform GRPO when many prompts are initially very difficult. On the same lines,  
\cite{tajwar2026maximum} connects this maximum-likelihood viewpoint to finite 
sampling by deriving rollout-indexed surrogate objectives that interpolate 
between a linear and logarithmic surrogate as the number of sampled responses 
increases.

\paragraph{Future directions:}
RLVR raises several natural questions for high-dimensional statistics. Recent theoretical work 
has begun to derive sample-complexity guarantees and identify base-model-dependent statistical 
barriers for policy-gradient post-training in stylized autoregressive settings 
\cite{mousavi2026post}. Extending such analyses to high-dimensional estimation 
regimes, modern RLVR surrogates, and group-based policy-gradient algorithms remains largely open.

\begin{itemize}
    \item \textbf{Surrogate design and statistical efficiency.} Can one develop tractable 
    generative models for prompts, answers, and responses under which the sample complexity, 
    generalization error, and optimality of different surrogate functions $F$ can be analyzed 
    as a function of the prompt distribution, model class, and number of sampled responses? 

    \item \textbf{RLVR versus supervised fine-tuning.} RLVR uses weak but scalable verifiable 
    feedback, whereas supervised fine-tuning uses richer labels. A basic open question is to 
    characterize when binary verifiable rewards are statistically sufficient, when richer 
    supervision is fundamentally necessary, and how the two forms of training could be combined.

\item \textbf{Noisy and imperfect feedback.} In practice, the verifier may itself be noisy or biased, for example when 
using an LLM-as-a-judge. While noisy or imperfect verification has begun to be studied  \cite{cai2025reinforcement,xu2025tinyv}, it remains open to characterize how 
false positives and false negatives affect the bias and variance of policy-gradient estimators, 
and whether the resulting dynamics remain consistent with the intended surrogate objective or 
instead optimize a noise-distorted one. Such an understanding could inform the design of robust 
RLVR algorithms and clarify when verifier imperfections lead to reward hacking \cite{amodei2016concrete}.
\end{itemize}

\section{Discussion}\label{ssec:discussion}

We have presented a selection of key research directions and open questions in high-dimensional statistics. Because the field is broad and rapidly evolving, many important areas have not received the attention they deserve. We now highlight a few notable additional directions and provide references for readers interested in exploring them further.

\subsection{False discovery rate control}

Controlling the false discovery rate (FDR) in high-dimensional settings is generally not a natural consequence of many existing model selection procedures \cite{su2017false}. The knockoff framework \cite{barber2015controlling} addresses this difficulty by augmenting the original variables with carefully constructed ``fake'' variables that act as negative controls, thereby enabling FDR control on top of a wide range of existing methods. Since its introduction, a substantial literature has extended and refined knockoffs for high-dimensional regression \cite{candes2018panning,dai2016knockoff,BatesEtAl2019Metro,barber2020robust,FanGaoLv2024ARK,RomanoEtAl2020DeepKnockoffs,katsevich2019multilayer}. However, a number of open questions remain. The first is generating knockoffs that yield high power. While existing approaches rely on MCMC \cite{BatesEtAl2019Metro} or generative models \cite{RomanoEtAl2020DeepKnockoffs}, further development is required to effectively handle complex high-dimensional data. Another important direction concerns the double robustness of knockoffs \cite{FanGaoLv2024ARK,FanGaoLvXu2025MomentsMatching}, which investigates whether FDR control holds when the distribution of the feature vector $\bm{X}$ is unknown or misspecified, provided the conditional model $\bm{Y} | \bm{X}$ is accurately specified. Finally, improving the power of knockoff-based inference continues to be an active area of research \cite{weinstein2023power,li2021whiteout,luo2025improving}. An ongoing challenge is to optimally combine knockoffs with other procedures to achieve strong power performance.

\subsection{Inference and uncertainty quantification}
While high-dimensional statistics initially focused on estimation and prediction under structural assumptions such as sparsity, a parallel and now substantial literature has developed surrounding statistical \emph{inference} in high-dimensional settings. In such settings, classical inferential procedures break down: regularization methods often introduce non-negligible bias, asymptotic normality of naive estimators fails, and model selection can invalidate inference which ignores model selection. A deeper understanding of these issues has led to the development of \emph{debiasing} methods for inference on low-dimensional parameters in the high-dimensional setting \cite{vandeGeer2014,ZhangZhang2014}. These methods are deeply connected to classical work on semi-parametric statistics~\cite{BickelKlaassenRitovWellner1993,ChernozhukovDML}. An alternative line of work has studied \emph{post-selection} inference, seeking valid inference conditional on a data-driven model selection step \cite{lee2016}. These inferential methods often require stringent assumptions for their validity, which in turn has led to a study of fundamental limits~\cite{CaiGuo2017} and the development of various ``assumption-light methods''~\cite{Rinaldo2019,hulc,universal} that are valid under weaker assumptions but often conservative in practice. Bridging the gaps between methods that are provably optimal under strong assumptions and those that are provably valid under weak assumptions is an active and exciting frontier of research.

\subsection{Tensors}\label{ssec:tensor}
%\textcolor{red}{Anru, could you please write a few paragaphs on this topic?}
High-dimensional tensor (or high-order) data are observations naturally represented as large multi-way arrays with three or more modes, and are increasingly prevalent in biology, medicine, psychology, education, computational imaging, and machine learning. Such data differ qualitatively from vectors and matrices: their scientific meaning is often encoded in multi-linear structure across modes (e.g., subject $\times$ feature $\times$ time), so methods and theory developed for matrices do not automatically extend. A common workaround is to reduce tensors to matrices or vectors via matricization or vectorization, but this can obscure intrinsic multi-way patterns and lead to statistically suboptimal analyses. On the computational side, tensors also introduce new barriers: many basic matrix notions do not admit tractable high-order analogues, and naive extensions of operator norms, singular values, and eigenvalues are NP-hard to compute in general \cite{hillar2013most}. These considerations motivate tensor-specific statistical methodology that is closely tied to the computational-statistical trade-off discussed in Section \ref{sec:trade-off}. Representative problems in high-dimensional tensor analysis include tensor PCA/SVD and related dimension reduction \cite{richard2014statistical,zhang2018tensor}, tensor completion \cite{montanari2018spectral,xia2021statistically}, tensor clustering \cite{sun2019dynamic}, and tensor regression \cite{zhou2013tensor}. Tensor methods also underpin broader inferential tasks such as interaction pursuit \cite{hao2014interaction}, independent component analysis \cite{auddy2025large}, learning latent variable models \cite{anandkumar2014tensor}, mixture model estimation \cite{doss2023optimal}, and Markov process estimation \cite{zhou2022optimal}. These ideas have enabled analyses across diverse applications, including single-cell high-throughput chromatin conformation capture (scHi-C) data \cite{park2024joint}, genome-wide association studies \cite{hore2016tensor}, tensor time series \cite{chen2022factor}, computational imaging \cite{han2020optimal}, and longitudinal sequencing data \cite{ma2023tensor,shi2024tempted}. Readers are referred to recent surveys and a book for broader introductions and overviews \cite{sun2021tensors,bi2021tensors,auddy2025tensors,KoldaBaderTensorBook}. Looking forward, major challenges include developing methods robust to misspecified structures, enabling scalable inference and uncertainty quantification, and integrating tensor methods with modern AI in ways that improve both performance and interpretability.

\subsection{Posterior contraction in high dimensions}

In the Bayesian framework, a statistical model $\bm{X} \sim \mathcal{P}_{\bm{\theta}}$ is equipped with a prior $\bm{\theta} \sim \Pi$, and posterior contraction refers to the property that
$
\Pi\big(L(\bm{\theta},\bm{\theta}^*) > \epsilon \mid \bm{X}\big) \to 0
$
in $P_{\bm{\theta}^*}$-probability, where $L(\cdot, \cdot)$ denotes a discrepancy measure. The primary tool for establishing contraction rates is the prior mass and testing framework \cite{barron1999consistency,ghosal2000convergence,shen2001rates,ghosal2017fundamentals}, which is designed for intrinsic loss functions such as the Hellinger distance \cite{ghosal2000convergence,GineNickl2011LrPosterior}. Over the past two decades, this theory has been extended to high-dimensional and nonparametric models \cite{castillo2012needles,castillo2015bayesian,gao2020general,banerjee2021bayesian,castillo2024bayesian}, demonstrating that suitably designed priors yield optimal contraction rates and often valid frequentist coverage for credible sets \cite{szabo2015frequentist}. However, the existing framework remains largely limited to intrinsic loss functions, and extending posterior contraction theory to non-intrinsic losses remains a major open challenge.

This difficulty is particularly pronounced in semiparametric and functional estimation problems, where the object of interest is a low-dimensional functional embedded in a high-dimensional or nonparametric model. Standard contraction theory typically yields rates for the joint model under Hellinger loss, but not for directional losses on the functional of interest. For example, in canonical correlation analysis, posterior contraction for the joint covariance matrix does not directly translate into rates for the leading canonical vectors.  In robust estimation under Huber contamination, specifying a prior on both the parameter and the contamination distribution is generally infeasible, even though simple frequentist estimators can achieve optimal rates \cite{chen2018robust}. Similarly, in a regression model, optimal inference for the noise level or a given regression coefficient is considerably more delicate than global estimation of the entire regression function \cite{yang2019posterior}. Related challenges arise in non-intrinsic losses such as the supremum norm in nonparametric regression and density estimation \cite{HoffmannRousseauSchmidtHieber2015AdaptivePosterior,YooGhosal2016SupNorm} or the operator norm in covariance estimation \cite{gao2016bernstein}, where estimation can be reinterpreted as simultaneous recovery of a collection of linear functionals \cite{Castillo2014SupNorm}. Finally, prior construction in high-dimensional settings often depends critically on the specific inferential target, leading to different priors for different functionals—a theoretical necessity in some cases but a practical obstacle in applied Bayesian analysis. Despite some recent progress based on posterior correction \cite{yiu2025semiparametric}, bridging this gap between general posterior contraction theory, functional estimation, and practical prior design remains an important open direction.

\subsection*{Acknowledgments} 
First and foremost, we would like to thank the editors of Statistical Science and Professor Robert Tibshirani for their initiative in organizing a collection of articles on open problems in statistics. Their initiative brought the authors of this paper together.

Arian Maleki would like to thank the National Science Foundation for support under Grants DMS-2515716 and DMS-2210506. Chao Gao is grateful for support from NSF Grants ECCS-2216912 and DMS-2310769, and an Alfred Sloan fellowship. Jason M. Klusowski is grateful for support from the National Science Foundation through
NSF CAREER DMS-2239448 and the Alfred P. Sloan Foundation through a Sloan Research Fellowship. Ali Shojaie gratefully acknowledges support from the National Institutes of Health through grant RF1-AG090462. Weijie Su was supported in part by NSF grant DMS-2310679 and would like to thank Zhimei Ren and Yiqiao Zhong for helpful discussions. Anru Zhang is grateful for support from the NIH Grants R01HL169347 and R01HL168940. Subhabrata Sen gratefully acknowledges support from NSF (DMS CAREER 2239234), ONR (N00014-23-1-2489) and AFOSR (FA9950-23-1-0429), and thanks Rahul Mazumder for helpful discussions. Verena Zuber gratefully acknowledges the United Kingdom Research and Innovation Medical Research Council grant MR/W029790/1. Christos Thrampoulidis acknowledges support from NSERC.

\end{sloppypar}

\bibliographystyle{alpha}
\bibliography{sample}
%\bibliography{references-trade-off}

\end{document}